\documentclass[english]{extarticle}
\usepackage[T1]{fontenc}
\usepackage[latin9]{inputenc}
\usepackage{color}
\usepackage{units}
\usepackage{textcomp}
\usepackage{amsmath}
\usepackage{amssymb}
\usepackage{graphicx}
\usepackage{esint}
\usepackage[authoryear]{natbib}

\makeatletter

\newcommand{\lyxmathsym}[1]{\ifmmode\begingroup\def\b@ld{bold}
  \text{\ifx\math@version\b@ld\bfseries\fi#1}\endgroup\else#1\fi}

\numberwithin{equation}{section}
\numberwithin{figure}{section}

\usepackage{babel}

\makeatother

\usepackage{babel}
\begin{document}

\title{\textcolor{black}{The coupon collector urn model with unequal probabilities
in ecology and evolution}}

\author{Zoroa, N.$^{1}$, Lesigne, E.$^{2}$, Fernández-Sáez, M.J.$^{1}$,
Zoroa, P.$^{1}$ \& Casas, J.$^{3}$}

\maketitle
$^{1}$Departamento de Estadística e Investigación Operativa, Facultad
de Matematicas, Universidad de Murcia, 30071, Murcia, SPAIN. $^{2}$Université
François-Rabelais, CNRS, LMPT UMR7350, Tours, FRANCE. $^{3}$Université
de Tours \& Institut Universitaire de France Institut de Recherche
en Biologie de l'Insecte, IRBI UMR CNRS 7261 37200, Tours, FRANCE
\begin{abstract}
The sequential sampling of populations with unequal probabilities
and with replacement in a closed population is a recurrent problem
in ecology and evolution. Many of these questions can be reformulated
as urn problems, often as special cases of the coupon collector problem,
most simply expressed as the number of coupons that must be collected
to have a complete set. We aimed to apply the coupon collector model
in a comprehensive manner to one example - hosts (balls) being searched
(draws) and parasitized (ball color change) by parasitic wasps - to
evaluate the influence of differences in sampling probabilities between
items on collection speed. 

Based on the model of a complete multinomial process over time, we
define the distribution, distribution function, expectation and variance
of the number of hosts parasitized after a given time, as well as
the inverse problem, estimating the sampling effort. We develop the
relationship between the risk distribution on the set of hosts and
the speed of parasitization and propose a more elegant proof of the
weak stochastic dominance among speed of parasitization, using the
concept of Schur convexity and the \textquotedblleft{} Robin Hood
transfer\textquotedblright{} numerical operation.

Numerical examples are provided and a conjecture about strong dominance
is proposed. The speed at which new items are discovered is a function
of the entire shape of the sampling probability distribution. The
sole comparison of values of variances is not sufficient to compare
speeds associated to two different distributions, as generally assumed
in ecological studies. 
\end{abstract}

\paragraph*{2010 Mathematics Subjet Classification:$\;$60C05, 60E15, 92B05}

\section{Introduction}

The description of sequential sampling of a population of individuals
for which the probability of being selected does not vary until a
specific event, such as the collection of all or some types of individuals
or a specific subgroup of the population, occurs is a common problem
in ecology and evolution studies. In probability theory, such problems
are often treated as urn problems, generally as the \textquotedblleft coupon
collector problem\textquotedblright{} (CCP). The CCP is a mathematical
model that belongs to the family of urn problems that can be formulated
as follows: A company issues coupons of different types, each with
a particular probability of being issued. The object of interest is
the number of coupons that must be collected to obtain a full collection.
This problem has been widely studied. The first findings concerned
the classical problem in which all coupons are equally likely to be
obtained \citep{Feller}. Rapid advances have been made in this field
\citep{BonehHofri,FlajoletBirthdayParadox,Anceaume2015}, but they
have gone largely unnoticed by most scientists working in ecology
and evolution. This is partly due to difficulties in making the correct
analogies, partly due to a lack of worked examples and partly because
each field devises its own vocabulary, procedures and formalism. In
ecological sciences, for example, a vibrant field of theoretical and
applied ecological statistics developed in the 1950s from the repeated
sampling of populations to estimate biodiversity richness \citep{McArthur,Simpson}.
This field could greatly benefit from the latest advances in the CCP
\citep{BungeFitzpatrick1993,HuilletParoissin2013}. Related problems
deal with the relative abundance of species from a community containing
many species \citep{Dennehy2009}, or the sampling effort required
to achieve a particular level of coverage \citep{NealMoriary2009}.
Increases in the number of new hosts being infected or superinfected
are a topic of great importance in population dynamics and epidemiology
\citep{DaleyGaniGani2001,Lloyd-SmithSchreiberKoppGetz(2005),keelingRohani2008}.
Several of the questions posed in capture-recapture studies relate
to the coupon collector problem. Occupancy problems and related capture-recapture
techniques are, indeed, defined as problems in which the probability
of a given species occupying a given state at a given time must be
determined (see the review \citealp{Baileyetal.2013} and the paper
\citealp{Hernandez-SuarezHiebeler2011}). In ethological sciences,
the estimation of a repertoire of signals in animal communication
is considered as a form of the CCP, because vocal repertoire size
is a key behavioral indicator of the complexity of the vocal communication
system in birds and mammals \citep{Kershenbaumetal.2015}. In genetics
and evolution, the coupon collector problem has been recognized as
such only occasionally, despite these fields having generated some
of the most elegant theorems and uses of other urn processes \citep{Ewens1972,Donelly1986,Doumas 2015}.
Indeed, the coupon collector problem has been used in the context
of exhaustive haplotype sampling in phylogeography, \citep{Dixon2006},
determining the number of beneficial mutations as a function of sequence
lines \citep{Tenaillonetal.2012}, and estimations of the size of
the library required to target a particular percentage of the non-essential
genome displaying a given property \citep{Vandewalleetal.2015}, for
example.

Urn models have been much more widely used for modelling host-parasitoid
systems than in other topics of ecology. We therefore used the biological
context and formalism of parasitism by parasitic wasps, as the results
obtained with this system can easily be extended to other ecological
and evolutionary contexts. Parasitic wasps search for insects hosts,
such as caterpillars, in which they lay a single, or multiple eggs.
In solitary wasp species, only a single wasp develops fully in a given
host. Parasitism can thus be formalized as a probabilistic dynamic
process with hosts as \textquoteleft balls\textquoteright{} and parasitoids
changing their \textquoteleft color\textquoteright{} by parasitizing
them. In work beginning more than a century ago \citep{Fiske1910,Thompson1924},
the pioneering population dynamicists assumed that hosts were found
and attacked on successive occasions governed by exponential laws
in continuous time. The number of draws was thus considered to be
random and the number of eggs for a given host was assumed to follow
a Poisson law \citep{Montovan}. If we assume that the number of draws
is fixed, then the distribution of the number of eggs for a given
host is binomial, but closely approximates a Poisson distribution
in large host populations. The proportion of the population without
eggs (the zero class) is of particular interest, because these hosts
survive parasitism and produce offspring for the next generation.
In field studies however, observed distributions are generally more
aggregated than would be expected under the assumed Poisson distribution
\citep{Hemericketal.2002}. Aggregation is interpreted as the result
of heterogeneity in the risk of being found, due to differences in
location, accessibility, appearance, color, developmental stage or
any other trait \citep{Hassell2000,MurdochBriggsNisbet2013}. The
risk distribution greatly influences the stability of the host-parasitoid
system and has been widely studied \citep{May1978,IvesSchoolerJagarKnutesonGrbicSettle1999,SinghMurdochNisbet2009}. 

All these models focus on the distribution of eggs over the entire
population of hosts, after a given time or a given number of draws
(Figure1). However, the use of this distribution greatly decreases
the amount of information available, as it collapses individual host
histories. Parasitism is a multinomial process (Figure1),
\begin{figure}
\includegraphics[width=12cm]{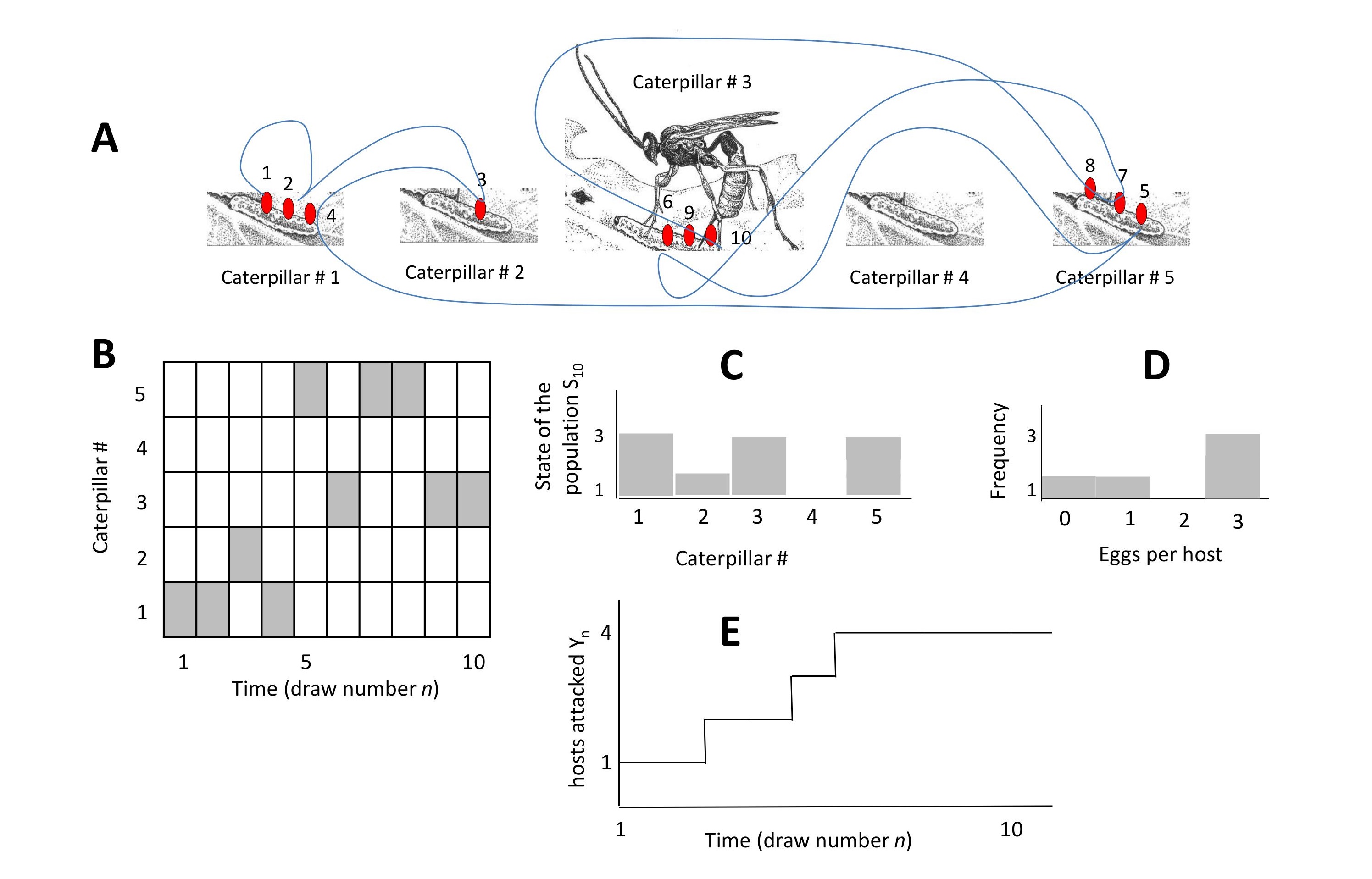}

Figure 1. Attacks of caterpillar larvae hosts by parasitic wasps as
an urn process in discrete time $n$. A wasp oviposit 10 eggs among
five hosts ($n=10$) (A). The outcome of the fundamental multinomial
process (B) is summarized in the marginal distribution of the number
of eggs per individual host at a given time $S_{10}=(3,1,3,0,3)$
(C), in the frequency distribution of eggs among hosts (D) and in
the number of hosts attacked over time $Y_{n}$ (E) 
\end{figure}
 in which time corresponds to host draws. Its dynamics determines,
for example, the percentage of hosts parasitized at the end of the
season, the opportunity and time at which alternative pest control
methods need to be deployed in supplement in biological control with
parasitoids, and the time required to achieve a given degree of control
by parasitic wasps. In the present paper, we aimed to model parasitism
as a multinomial urn process over time and we study the speed of parasitization
(Figure1). We consider host encounters followed by oviposition without
discrimination. The parasitism process described above can be considered
as a coupon collector problem. In this case, there is a finite population
of hosts differing in appearance, location, developmental stage or
other factors. This heterogeneity results in different probabilities
of hosts being found by parasitoids. These probabilities, $p_{h}$
for host $h$, do not change over time. Our work therefore entails
describing in depth the coupon collector problem, highlighting unnoticed
analogies among previous works within the probability literature,
and comparing the influence of the degree of heterogeneity among hosts
on the speed of infection. We give a compact and hopefully more elegant
proof than previously known of the following fact : the more the distribution
$p$ on the set of hosts is heterogeneous, the more the (random) number
$Y$ of parasitized hosts after a given number of draws is small;
in other terms, there is a monotonic relationship between the majorization
relation on the set of probability distributions $p$ with the stochastic
dominance on the set of random numbers $Y$.

This paper is structured as follows. In Section 2, we define a succession
$S_{n}$, $n=1,2,...,$ of $N$-dimensional random variables describing
the state of the host population over time, in which time, $n,$ is
given by the number of attacks on the set of hosts. Each marginal
distribution of $S_{n}$ provides us information about a subset of
hosts, including, in particular, the $h$-th component representing
the number of times that host $h$ has been attacked by a parasitoid
between times 1 and $n$. In Section 3 we define the random variables
$Y_{n}$, $n=1,2,...,$ representing the number of parasitized hosts
after $n$ draws. We also compute the distribution, the distribution
function, the expectation and the variance of $Y_{n}$. We found no
examples of calculations of this value in previous studies and therefore
believe this aspect to be novel. We obtain the expected number of
draws required for all hosts in a given subset to be parasitized and
provide upper and lower bounds for this value in Section 4. In Section
5, we apply the results developed in previous sections to two particular
risk distributions on the set of hosts. We first use the uniform distribution,
and then a distribution corresponding to a host population with two
different kinds of hosts. We calculate the most relevant values for
each of these cases. In Section 6, we develop the relationship between
the speed of parasitization and the risk distribution in the set of
hosts. A narrower risk distribution is associated with faster parasitization.
Thus, parasitization is fastest when the risk distribution is uniform.
We highlight this finding with numerical examples in Section 7 and
propose a conjecture on strong stochastic dominance in Section 8.

\section{Modelling parasitism as an urn process}

We assume a finite population of hosts, constant for the entire duration
of the experiment. The parasitoid population is irrelevant, but we
assume that the number of eggs that can be laid in the host population
is not limiting. The situation is developed in successive stages or
draws. At each stage, a parasitoid attacks a host, in which it lays
an egg. \textit{\emph{The model is based on the fundamental assumption
that successive draws are mutually independent. The hosts differ in
appearance due to intrinsic qualities, and these differences modify
their probability of being attacked by a parasitoid. If the hosts
are named 1, 2, 3, \dots , $N$, then host $h$ has a probability
$p_{h}\geq0$ ($\sum p_{h}=1$) of being attacked by a parasitoid
in a draw. This probability does not change during the process. We
will say that $p_{1}$,$p_{2}$,...,$p_{N}$ or $(p_{1},p_{2},...,p_{N})$
is the risk distribution for the set of hosts $H=\left\{ 1,2,...,N\right\} $. }}

The underlying probability space of our model is $(\Omega,\mathcal{S},P)$,
where the elements of $\Omega$ are all the possible histories of
parasitism, that is $\Omega=H^{\mathbb{N}}$, equipped with its product
$\sigma$-algebra $\mathcal{S}$ and the probability $P$ given by
Kolmogorov\textasciiacute s Theorem: if we therefore fix $i_{1},i_{2},i_{3},\ldots,i_{n}$
in $H$, the probability of the event $\{\omega=i_{1}i_{2}i_{3}\ldots i_{n}j_{n+1}j_{n+2}\dots:\text{ for some \ensuremath{j_{k}}in \ensuremath{H},\ensuremath{\: k>n}}\}$
is $p_{i_{1}}p_{i_{2}}\ldots p_{i_{n}}$. When necessary, the vector
$(p_{1},p_{2},...,p_{N})$ will be denoted by a single letter $p$,
and the probability $P$ will be denoted by $P_{p}$.

We can describe this situation by defining a succession of random
variables, 
\begin{equation}
S_{n}=(S_{n1},S_{n2},...,S_{nN}),\quad n=1,2,...\label{eq:Sn}
\end{equation}
where $S_{nj}$ denotes the number of eggs in host $j$ after $n$
draws.

Variable $S_{n}$ represents the state of the host population after
$n$ draws, that is, the distribution of eggs over the total population
of hosts. If host $i$ was visited $r_{i}$ times between stages 1
and $n$, for $i=1,2,...,N$, then $S_{n}$ takes the value $(r_{1},r_{2},...,r_{N})$.
This variable has a multinomial distribution with parameters $n$,
$p_{1}$,$p_{2}$,...,$p_{N}$, that is, for every integers $r_{1},$
$r_{2}$,..., $r_{N}$, $0\leq r_{i}\leq n$, $r_{1}+r_{2}+...+r_{N}=n$,
\begin{equation}
P(S_{n1}=r_{1},S_{n2}=r_{2},...,S_{nN}=r_{N})=\frac{n!}{r_{1}!r_{2}!...r_{N}!}p_{1}^{r_{1}}p_{2}^{r_{2}}...p_{N}^{r_{N}}.\label{eq:Prob Sn1=00003D00003Dr1, Sn2=00003D00003Dr2, etc}
\end{equation}

The marginal distribution of $S_{n}$, for $i_{1},i_{2},...,i_{h}$
distinct elements of $\left\{ 1,2,...,N\right\} $ is given by 
\[
P(S_{ni_{1}}=r_{i_{1}},S_{ni_{2}}=r_{i_{2}},...,S_{ni_{h}}=r_{i_{h}})=
\]
\begin{equation}
\frac{n!}{r_{i_{1}}!r_{i_{2}}!...r_{i_{h}}!(n-\sum r_{i_{j}})!}p_{i_{1}}^{r_{i_{1}}}p_{i_{2}}^{r_{i_{2}}}...p_{i_{h}}^{r_{i_{h}}}q_{i_{1}i_{2}...i_{h}}^{n-\sum r_{i_{j}}},\label{eq:prob marginal}
\end{equation}
\[
\mathrm{with}\quad0\leq r_{i_{j}}\leq n,\; j=1,2,...,h,\; r_{i_{1}}+r_{i_{2}}+...+r_{i_{h}}\leq n,\; q_{i_{1}i_{2}...i_{h}}=1-\sum_{j=1}^{h}p_{i_{j}}{\color{red}.}
\]
This is the probability that, after $n$ draws host $i_{1}$ has been
visited $r_{i_{1}}$ times by the parasitoids, host $i_{2}$ $r_{i_{2}}$
times and host $i_{h}$ $r_{i_{h}}$ times, without considering the
rest of the hosts.

In particular, the component $S_{nh}$ of $S_{n}$ has a binomial
distribution with parameters $n$, $p_{h}$, 
\begin{equation}
P(S_{nh}=r)=\frac{n!}{r!(n-r)!}p_{h}^{r}q_{h}^{n-r},r=0,1,2,...,n\;\label{eq:prob component}
\end{equation}
\[
\mathrm{where}\; q_{h}=1-p_{h}=\sum_{i\neq h}p_{i}.
\]
This variable represents the state of host $h$ after $n$ draws.
Thus, $P(S_{nh}=r)$ is the probability that host $h$ has been attaked
$r$ times during the $n$ draws.

The expected value and variance of this random variable are 
\[
E(S_{nh})=np_{h},\quad Var(S_{nh})=np_{h}(1-p_{h}).
\]

Let $(e_{1},e_{2},...,e_{N})$ denote the canonical base of the space
$\mathbb{R}^{N}$. We emphasize that the process $(S_{n})_{n\geq1}$
is the random walk on $Z_{+}^{N}$ with independent increments obeying
the following law: $S_{n+1}-S_{n}=e_{k}$ with probability $p_{k}$.
The statistical behavior of this process is also very well known.

Note that, in this model, the sequence of random subsets of $H$,
describing the set of parasitized hosts over time, is a Markov chain,
and it is not difficult to give a precise description of its probability
transitions. However, it is not straightforward to study this Markov
chain directly.

\section{Number of parasitized hosts after $n$ draws}

Let $Y_{n}$ be the random variable representing the number of parasitized
hosts after $n$ draws, that is $Y_{n}=k$ if there are exactly $k$
parasitized hosts after $n$ draws. In this section we study this
random variable obtaining expressions for: its probability mass function
(\ref{eq:Pnk notacion conjuntos}), distribution function (\ref{eq:FuncDistribYn}),
expectation (\ref{eq:E(Yn)}) and variance (\ref{eq:VarYn}).

From now on, for any integer $h>0$ and real $x$, we write
\[
x^{(h)}=x(x-1)(x-2)...(x-h+1),\hspace*{1em}x^{(0)}=1,
\]

\begin{equation}
\dbinom{x}{h}=\dfrac{x^{(h)}}{h!}=\dfrac{x(x-1)...(x-h+1)}{h!},\;\mathrm{and}\;\dbinom{x}{0}=1.\label{eq:NumeroCombinatorio}
\end{equation}

The distribution and the distribution function of $Y_{n}$ have been
obtained in previous studies, see \citet{Anceaume2015}\textcolor{black}{{}
}.

By applying the general inclusion and exclusion principle, we find
that, for any distinct elements $j_{1},j_{2},...,j_{k}$ of $H$,
and denoting $p_{j_{1}}{}_{j_{2}...}{}_{j_{k}}=p_{j_{1}}+p_{j_{2}}+...+p_{j_{k}}$,
\[
P(\mathrm{the\: set\: of\: parasitized\: hosts\: after\:}n\:\mathrm{draws\: is\,}\left\{ j_{1},j_{2},...,j_{k}\right\} )=
\]
\[
p_{j_{1}}^{n}{}_{j_{2}...}{}_{j_{k}}-\underset{\left\{ i_{1},i_{2},...,i_{k-1}\right\} \subset\left\{ j_{1},j_{2},...,j_{k}\right\} }{\sum}p_{i_{1}}^{n}{}_{i_{2}...}{}_{i_{k-1}}+
\]
\[
\underset{\left\{ i_{1},i_{2},...,i_{k-2}\right\} \subset\left\{ j_{1},j_{2},...,j_{k}\right\} }{\sum}p_{i_{1}}^{n}{}_{i_{2}...}{}_{i_{k-2}}-...+(-1)^{k-1}\underset{i\in\left\{ j_{1},j_{2},...,j_{k}\right\} }{\sum}p_{i}^{n},
\]
from which we deduce that

\begin{center}
$P(Y_{n}=k)=\underset{\left\{ j_{1},j_{2},...,j_{k}\right\} \subset\left\{ 1,2,...,N\right\} }{\sum}p_{j_{1}}^{n}{}_{j_{2}...}{}_{j_{k}}-\left(\begin{array}{c}
N-k+1\\
N-k
\end{array}\right)\underset{\left\{ j_{1},j_{2},...,j_{k-1}\right\} \subset\left\{ 1,2,...,N\right\} }{\sum}p_{j_{1}}^{n}{}_{j_{2}...}{}_{j_{k-1}}$ 
\par\end{center}

\begin{center}
$+\left(\begin{array}{c}
N-k+2\\
N-k
\end{array}\right)\underset{\left\{ j_{1},j_{2},...,j_{k-2}\right\} \subset\left\{ 1,2,...,N\right\} }{\sum}p_{j_{1}}^{n}{}_{j_{2}...}{}_{j_{k-2}}-...$ 
\par\end{center}

\begin{center}
\begin{equation}
+(-1)^{k-1}\left(\begin{array}{c}
N-1\\
N-k
\end{array}\right)\underset{\left\{ j\right\} \subset\left\{ 1,2,...,N\right\} }{\sum}p_{j}^{n}.\label{eq:Pnk}
\end{equation}

\par\end{center}

Using the notation $p_{A}=\sum_{i\in A}\, p_{i}$ for any $A\subset H$,
this can be written in a more compact form 
\begin{equation}
P(Y_{n}=k)=\sum_{A\subset H,\left|A\right|\leq k}(-1)^{k-\left|A\right|}\dbinom{N-\left|A\right|}{k-\left|A\right|}\: p_{A}^{n},\quad\mathrm{for}\;0\leq k\leq N,\;\mathrm{and\;}k\leq n.\label{eq:Pnk notacion conjuntos}
\end{equation}
where $\left|A\right|$ denotes the number of elements of the set
$A$.

Let us now consider the distribution function of $Y_{n}$, 
\[
P(Y_{n}\leq k)=\sum_{j=1}^{k}P(Y_{n}=j)=\sum_{j=1}^{k}\sum_{A\subset H,\left|A\right|\leq j}(-1)^{j-\left|A\right|}\dbinom{N-\left|A\right|}{j-\left|A\right|}\: p_{A}^{n}=
\]
\begin{equation}
\sum_{A\subset H,\left|A\right|\leq k}\left(\sum_{i=0}^{k-\left|A\right|}(-1)^{i}\dbinom{N-\left|A\right|}{i}\right)\: p_{A}^{n}.\label{eq:FunciondistribucionYnk}
\end{equation}
As, for any integers $K$ and $k\geq0$ the equality 
\[
\sum_{i=0}^{k}(-1)^{i}\dbinom{K}{i}=(-1)^{k}\dbinom{K-1}{k}
\]
holds, we obtain 
\begin{equation}
P(Y_{n}\leq k)=\sum_{A\subset H,\left|A\right|\leq k}(-1)^{k-\left|A\right|}\dbinom{N-\left|A\right|-1}{k-\left|A\right|}\: p_{A}^{n},\quad k=1,2,...,N.\label{eq:FuncDistribYn}
\end{equation}

A similar expression can be seen in \citet{Anceaume2015} . From (\ref{eq:FuncDistribYn})
we can calculate the moments of $Y_{n}$. Let 
\[
m_{k}^{[n]}=\sum_{A\subset H,\left|A\right|=k}\: p_{A}^{n},
\]
for every $k\leq N$ 
\[
\sum_{j=1}^{k}P(Y_{n}\leq j)=\sum_{l=1}^{k}(-1)^{k-l}\dbinom{N-l-2}{k-l}\: m_{l}^{[n]},
\]
this gives, with $k=N-1$ and $k=N$ , 
\[
\sum_{j=1}^{N-1}P(Y_{n}\leq j)=m_{N-1}^{[n]},
\]
\[
\sum_{j=1}^{N}P(Y_{n}\leq j)=m_{N}^{[n]}+m_{N-1}^{[n]}=1+m_{N-1}^{[n]}.
\]
And, bearing in mind that 
\[
E(Y_{n})=\sum_{j=1}^{N}P(Y_{n}\geq j)=1+\sum_{j=1}^{N}P(Y_{n}>j)=1+N-\sum_{j=1}^{N}P(Y_{n}\leq j),
\]
we obtain the well known formula: 
\begin{equation}
E(Y_{n})=N-m_{N-1}^{\left[n\right]}=N-\sum_{i=1}^{N}(1-p_{i})^{n}.\label{eq:E(Yn)}
\end{equation}
We were unable to find any expression for $E(Y_{n}^{2})$ and the
variance of $Y_{n}$, in previous studies. These two quantities can
be obtained as follows. We compute, for $k\leq N$ 
\[
\sum_{t=1}^{k}\sum_{j=1}^{t}P(Y_{n}\leq j)=\sum_{l=1}^{k}(-1)^{k-l}\dbinom{N-l-3}{k-l}\: m_{l}^{[n]},
\]
then, for $k=N-2,\; N-1$ and $N$ we obtain 
\[
\sum_{t=1}^{N-2}\sum_{j=1}^{t}P(Y_{n}\leq j)=\: m_{N-2}^{[n]},
\]
\[
\sum_{t=1}^{N-1}\sum_{j=1}^{t}P(Y_{n}\leq j)=\: m_{N-2}^{[n]}+m_{N-1}^{[n]},
\]
and 
\[
\sum_{t=1}^{N}\sum_{j=1}^{t}P(Y_{n}\leq j)=\: m_{N-2}^{[n]}+2m_{N-1}^{[n]}+m_{N}^{[n]}.
\]
The last identity can be written: 
\[
\sum_{j=1}^{N}\dfrac{(N-j+1)(N-j+2)}{2}P(Y_{n}=j)=\: m_{N-2}^{[n]}+2m_{N-1}^{[n]}+m_{N}^{[n]},
\]
this gives 
\[
\dfrac{(N+1)(N+2)}{2}-\dfrac{2N+3}{2}E(Y_{n})+\dfrac{1}{2}E(Y_{n}^{2})=\: m_{N-2}^{[n]}+2m_{N-1}^{[n]}+m_{N}^{[n]},
\]
therefore 
\[
E(Y_{n}^{2})=2m_{N-2}^{[n]}-(2N-1)m_{N-1}^{[n]}+N^{2}=
\]
\[
2\sum_{1\leq i<j\leq N}(1-p_{i}-p_{j})^{n}-(2N-1)\sum_{i=1}^{N}(1-p_{i})^{n}+N^{2}
\]
and 
\begin{equation}
Var(Y_{n})=2\sum_{1\leq i<j\leq N}(1-p_{i}-p_{j})^{n}+\sum_{i=1}^{N}(1-p_{i})^{n}\left(1-\sum_{i=1}^{N}(1-p_{i})^{n}\right).\label{eq:VarYn}
\end{equation}

\section{\textcolor{black}{The number of draws required to reach a given level
of parasitism}}

The expected number of draws required for the parasitization of $k$
unparasitized hosts may be of considerable interest. For example,
we might want to know the expected number of draws required for $k$
of the hosts occupying a determined region, or with probabilities
of parasitization greater (or less) than a given value, etc., are
parasitized. We define below a random variable representing the number
of draws required for the event of interest to happen and we obtain
its expectation. We also describe the relationship between the random
variables defined here and the variables $Y_{n}$ defined in Section
3.

Let us consider that, at a given stage of the process, there is a
set $K\subset H$ of unparasitized hosts, this is our set of interest,
and the remaining hosts $H-K$ are or are not parasitized. Let us
use $X$ to denote the number of hosts in the set $H-K$ attacked
by the parasitoids before one of the hosts in $K$ is attacked.

As this process involves the repeating of independent trials, the
random variable $X$ follows a geometric distribution with parameter
$p=\sum_{i\in K}p_{i},$ (or a degenerate distribution if $K=H$).
Thus, 
\begin{equation}
E(X)=\frac{\sum_{i\in H-K}p_{i}}{\sum_{i\in K}p_{i}}.\label{eq:Expected draws before new parasitation}
\end{equation}

Now, let $k$ and $N_{\text{1}}$ be integers $1\leq k\leq N_{1}\leq N$.
Let $H_{1}$ be a subset of the set of hosts, $H$, and $H_{2}=H-H_{1}$,
$\left|H_{1}\right|=N_{1}$. We can assume that $H_{1}=\left\{ 1,2,...N_{1}\right\} $
without lost of generality.

If we consider the hosts of set $H_{1}$ to be unparasitized, then
we can define $T_{k,N_{1}}$ as the random number of draws required
to ensure that $k$ hosts of set $H_{1}$ have been parasitized. Its
expectation is the expected number of draws required for $k$ hosts
of set $H_{1}$ be parasitized. The case $H_{1}=H$ has been studied
before and different expressions for $E(T_{k,N})$ have been obtained.
We include these at the end of this section. In \citet{BonehHofri},
an expression is proposed for the particular case in which $k=N_{1}=N$.
The expression obtained here is more general.

Let $i_{1},i_{2},...,i_{k}$ be distinct elements of $H_{1}$. Let
$D_{i_{1}i_{2}...i_{k}}$ be the event defined by the fact that the
first $k$ hosts of set $H_{1}$ parasitized (i.e. attacked by a parasitoid
for first time) are hosts $i_{1},i_{2},...,i_{k}$ and are parasitized
in the precise order $i_{1},i_{2},...,i_{k}$. In other words, some
of the hosts of set $H_{2}$ may be attacked first, followed by host
$i_{1}$. Next, some hosts of $H_{2}\cup\left\{ i_{1}\right\} $ may
be attacked, followed by host $i_{2}$, etc. Let $p=\sum_{i\in H_{1}}p_{i}$,
$q=1-p=\sum_{i\in H_{2}}p_{i}$. Then

\[
P(D_{i_{1}i_{2}...i_{k}})=P(\mathrm{first\; host\; of\;\mathit{H\mathrm{_{1}}}\; parasitized\; is\;}i_{1})
\]
\[
P(\mathrm{second\; host\; of\;\mathit{H\mathrm{_{1}}}\; parasitized\; is\;}i_{2}\;|\;\mathrm{first\; host\; of\;\mathit{H\mathrm{_{1}}}\; parasitized\; was\;}i_{1})\;...
\]
\[
P(k-\mathrm{th}\hspace*{1em}\mathrm{host\, of\;\, H\mathrm{_{1}}\, parasitized\; is}\: i_{k}|\,\mathrm{first\;(\mathit{k}-1)\mathrm{\, hosts\, of\;\mathit{H\mathrm{_{1}}}\, parasitized}\, were\,}i_{1},i_{2},...,i_{k-1}).
\]
Both in the case $q=0$ ($H_{1}=\left\{ 1,2,...,N\right\} $) and
the case $q>0$

\[
P(\mathrm{first\; parasitized\; host\; of\;\mathit{H\mathrm{_{1}}}\; is\;}i_{1})=\dfrac{p_{i_{1}}}{1-q},
\]
where 
\[
q=\sum_{i\in H_{2}}p_{i}>0.
\]
For the rest of the factors 
\[
P(\mathrm{\mathit{h-\mathrm{th}}\, parasitized\, host\, of\,\mathit{H\mathrm{_{1}}}\, is}\, i_{h}|\,\mathrm{first\, parasitized\; hosts\; of\;\mathit{H\mathrm{_{1}}}\, were\,}i_{1},\; i_{2},...,i_{h-1})=
\]

\[
\sum_{r=0}^{\infty}p_{i_{h}}(q+\sum_{j=1}^{h-1}p_{i_{j}})^{r}=\dfrac{p_{i_{h}}}{1-q-\sum_{j=1}^{h-1}p_{i_{j}}},\; h=1,2,...,k,\;\mathrm{where}\; q=\sum_{i\in H_{2}}p_{i},
\]
therefore, 
\begin{equation}
P(D_{i_{1}i_{2}...i_{k}})=\dfrac{\prod_{j=1}^{k}p_{i_{j}}}{p(p-p_{i_{1}})(p-p_{i_{1}}-p_{i_{2}})...(p-\sum_{j=1}^{k-1}p_{i_{j}})}.\label{eq:P(Di1i2...in1)}
\end{equation}

$\:$

Let $\Pi_{k}$ be the set of all $k$-permutations of 1, 2, ..., $N_{1}$.
Then the events $D_{i_{1}i_{2}...i_{k}}$ with $(i_{1}i_{2}...i_{k})\in\Pi_{k}$
constitute a partition of $\Omega$, i.e. $D_{i_{1}i_{2}...i_{k}}\cap D_{j_{1}j{}_{2}...j_{k}}=\oslash$
if $i_{1}i_{2}...i_{N_{1}}\neq j_{1}j{}_{2}...j_{N_{1}}$ and 
\[
\sum_{(i_{1}i_{2}...i_{k})\in\Pi_{k}}P(D_{i_{1}i_{2}...i_{k}})=1.
\]
We can then write $E(T_{k,N_{1}})$ as follows,

\begin{equation}
E(T_{k,N_{1}})=\sum_{(i_{1}i_{2}...i_{k})\in\Pi_{k}}E(T_{k,N_{1}}|D_{i_{1}i_{2}...i_{k}})P(D_{i_{1}i_{2}...i_{k}}).\label{eq:E(TkN1)=00003D00003DSumE(TkN1/D) P(D)}
\end{equation}

To compute the conditional expectations $E(T_{k,N_{1}}|D_{i_{1}i_{2}...i_{k}})$,
let us denote by $X_{h}$ the random variable representing the number
of draws elapsed after $h-1$ hosts of the set $H_{1}$ being parasitized
and before a new host of the set $H_{1}$ is parasitized, $1\leq h\leq k$.
We can then write 
\begin{equation}
T_{k,N_{1}}=X{}_{1}+1+X_{2}+1+...+X_{k}+1=X{}_{1}+X_{2}+...+X_{k}+k.\label{eq:TkN1=00003D00003DSumXh}
\end{equation}
and therefore

\begin{equation}
E(T_{k,N_{1}}|D_{i_{1}i_{2}...i_{k}})=\sum_{h=1}^{k}E(X_{h}|D_{i_{1}i_{2}...i_{k}})+k\label{eq:E(T/D)=00003D00003DSum E(Xh/D)}
\end{equation}
but

\begin{equation}
E(X_{h}|D_{i_{1}i_{2}...i_{k}})=E\left(X_{h}|\mathrm{already\: parasitized\: hosts\: are\: those\: of\:}H_{2}\:\mathrm{and\:}i_{1}i_{2}...i_{h-1}\right)\label{eq:Nada}
\end{equation}
One direct application of (\ref{eq:Expected draws before new parasitation})
would then be: 
\begin{equation}
E(X_{h}|D_{i_{1}i_{2}...i_{k}})=\dfrac{q+\sum_{j=1}^{h-1}p_{i_{j}}}{p-\sum_{j=1}^{h-1}p_{i_{j}}}\quad\mathrm{for\;}h=1,2,...,k.\label{eq:E(Xh/D)}
\end{equation}

From (\ref{eq:E(T/D)=00003D00003DSum E(Xh/D)}) and (\ref{eq:E(Xh/D)})
\[
E(T_{k,N_{1}}|D_{i_{1}i_{2}...i_{k}})=\left(\sum_{h=1}^{k}\dfrac{q+\sum_{j=1}^{h-1}p_{i_{j}}}{p-\sum_{j=1}^{h-1}p_{i_{j}}}\right)+k=\sum_{h=1}^{k}(\dfrac{q+\sum_{j=1}^{h-1}p_{i_{j}}}{p-\sum_{j=1}^{h-1}p_{i_{j}}}+1)=
\]
\[
\frac{1}{p}+\frac{1}{p-p_{i_{1}}}+\frac{1}{p-p_{i_{1}}-p_{i_{2}}}+...+\dfrac{1}{p-\sum_{j=1}^{k-1}p_{i_{j}}}=
\]
\begin{equation}
\dfrac{1}{1-q}+\dfrac{1}{1-q-p_{i_{1}}}+...+\dfrac{1}{1-q-\sum_{j=1}^{k-1}p_{i_{j}}},\label{eq:E(Z|Di1 i2 ...)}
\end{equation}
where 
\[
q=\sum_{i\in H_{2}}p_{i}=\sum_{i=N_{1}+1}^{N}p_{i}\quad\mathrm{and}\quad p=\sum_{i\in H_{1}}p_{i}=\sum_{i=1}^{N_{1}}p_{i}.
\]

$\:$

Bearing in mind (\ref{eq:E(TkN1)=00003D00003DSumE(TkN1/D) P(D)}),
(\ref{eq:P(Di1i2...in1)}) and (\ref{eq:E(Z|Di1 i2 ...)}) we can
state the following:

$\:$

$\mathbf{Proposition\:4.1.}$ The expected value of $T_{k,N_{1}}$
is 
\begin{multline*}
E(T_{k,N_{1}})=\sum_{(i_{1}i_{2}...i_{k})\in\Pi_{k}}\left(\frac{1}{p}+\frac{1}{p-p_{i_{1}}}+\frac{1}{p-p_{i_{1}}-p_{i_{2}}}+...+\dfrac{1}{p-\sum_{j=1}^{k-1}p_{i_{j}}}\right)
\end{multline*}

\begin{multline}
\dfrac{\prod_{j=1}^{k}p_{i_{j}}}{p(p-p_{i_{1}})(p-\sum_{j=1}^{2}p_{i_{j}})...(p-\sum_{j=1}^{k-1}p_{i_{j}})},\label{eq:E(Zk)}
\end{multline}
where $\Pi_{k}$ is the set of all $k$-permutations of set $\left\{ 1,2,...,N_{1}\right\} $,
i.e. the arrangements of length $k$ of different elements of $\left\{ 1,2,...,N_{1}\right\} $.

$\:$

Thus, $E(T_{k,N_{1}})$ given by (\ref{eq:E(Zk)}) is the expected
number of draws required for $k$ hosts of a set of unparasitized
hosts $H_{1}\subset H$ with cardinality $N_{1}$, to be parasitized.
This value is generally difficult to obtain because the number of
terms required for its computation is the number of $k$-permutations
of 1, 2, ..., $N_{1}$, that is $N_{1}^{(k)}=N_{1}(N_{1}-1)...(N_{1}-k+1)$.
This value is huge when $N_{1}$ and $k$ are large. It is therefore
important to obtain upper and lower bounds for this value. 

$\vphantom{}$

$\mathbf{Proposition\:4.2.}$ Let $k$ be given and $p_{1}$, $p_{2}$,
...,$p_{N_{1}}$ be real numbers satisfying $p_{1}\geq p_{2}\geq...\geq p_{N_{1}}$.
Then, the maximum of $E(T_{k,N_{1}}|D_{i_{1}i_{2}...i_{k}})$ defined
by (\ref{eq:E(Z|Di1 i2 ...)}) over all possible choices of the $k$-subsets
$\left\{ i_{1},i_{2},...,i_{k}\right\} $ of $H_{1}$ is 
\begin{equation}
E(T_{k,N_{1}}|D_{1,2,...,k})=\dfrac{1}{\sum_{i=1}^{N_{1}}p_{i}}+\dfrac{1}{\sum_{i=2}^{N_{1}}p_{i}}+...+\dfrac{1}{\sum_{i=k}^{N_{1}}p_{i}},\label{eq:D12...N1}
\end{equation}
and the minimum is 
\begin{equation}
E(T_{k,N_{1}}|D_{N_{1},N_{1}-1,...,N_{1}-k+1})=\dfrac{1}{\sum_{i=1}^{N_{1}}p_{i}}+\dfrac{1}{\sum_{i=1}^{N_{1}-1}p_{i}}+...+\dfrac{1}{\sum_{i=1}^{N_{1}-k+1}p_{i}}.\label{eq:DN1...321}
\end{equation}

$\vphantom{}$

$\mathbf{Proof.}$ From hypothesis $p_{1}\geq p_{2}\geq...\geq p_{N_{1}}$,
it follows directly that 
\begin{equation}
\sum_{i=h}^{N_{1}}p_{i}\leq\sum_{j=h}^{N_{1}}p_{i_{j}}\leq\sum_{i=1}^{N_{1}-h+1}p_{i},\; h=1,2...,N_{1},\label{eq:Sum pi < Sum pij}
\end{equation}
then 
\[
E(T_{k,N_{1}}|D_{1,2,...,k})=\frac{1}{p}+\frac{1}{p-p_{1}}+\dfrac{1}{p-\sum_{i=1}^{2}p_{i}}+...+\dfrac{1}{p-\sum_{i=1}^{k-1}p_{i}}\geq
\]
\[
\frac{1}{p}+\frac{1}{p-p_{i_{1}}}+\frac{1}{p-p_{i_{1}}-p_{i_{2}}}+...+\dfrac{1}{p-\sum_{j=1}^{k-1}p_{i_{j}}}\geq
\]
\[
\frac{1}{p}+\frac{1}{p-p_{N_{1}}}+\dfrac{1}{p-\sum_{i=N_{1}-1}^{N_{1}}p_{i}}+...+\dfrac{1}{p-\sum_{i=N_{1}-k+2}^{N_{1}}p_{i}}=
\]
\[
E(T_{k,N_{1}}|D_{N_{1},N_{1}-1,...,N_{1}-k+1})
\]
and the proof is complete.

$\vphantom{}$

$\mathbf{Proposition\:4.3.}$ Let $p_{1}$, $p_{2}$, ..., $p_{N_{1}}$
be real numbers satisfying $0\leq p_{i}\leq1$, for $i=1,2,...,N_{1}$
and $p_{1}\geq p_{2}\geq...\geq p_{N_{1}}$. It is then true that
\[
E(T_{k,N_{1}}|D_{1,2,...,k})\geq E(T_{k,N_{1}})\geq E(T_{k,N_{1}}|D_{N_{1},N_{1}-1,...N_{1}-k+1})
\]
In other words, $E(T_{k,N_{1}}|D_{1,2,...,k})$ and $E(T_{k,N_{1}}|D_{N_{1},N_{1}-1,...N_{1}-k+1})$
are upper and lower bounds, respectively, for the expected number
of draws required for $k$ hosts of the set $H_{1}$ to be parasitized.

Furthermore, the mode of the distribution on the events $D_{i_{1}i_{2}...i_{k}},\;(i_{1},i_{2},...,i_{k})\in\Pi_{k}$,
is $D_{1,2,...,k}$, i.e. the order of parasitism of $k$ hosts mostly
likely to occur is $1,2,...,k.$

$\vphantom{}$

$\mathbf{Proof.}$ The first part of this proposition is a straightforward
consequence of the previous proposition.

The second part comes directly from the fact that 
\[
P(D_{1,2,...,k})\geq P(D_{i_{1}i_{2}...i_{k}})\;\mathrm{for}\;(i_{1}i_{2}...i_{k})\in\Pi_{k}.
\]
which follows from (\ref{eq:P(Di1i2...in1)}) and (\ref{eq:Sum pi < Sum pij}). 

$\:$

Propositions 4.2 and 4.3 prove that, if $p_{1}\geq p_{2}\geq...\geq p_{N_{1}}$,
then the most likely order of parasitization of $k$ hosts in $H_{1}$
is the preferential order 1, 2, ..., $k$. Moreover the shortest scenario
(in terms of expectation) for the parasitization of $k$ hosts of
$H_{1}$ is the sequence extending from the least likely host, $N_{1}$,
to the most likely host, $N_{1}-k+1$, in the correct order. The longest
scenario (in terms of expectation) for the parasitization of $k$
hosts of $H_{1}$ extends from the most likely, 1, to the least likely
host, $k$, in the correct order.

These results can be intuitively explained as follows; let us suppose
that host 1 is parasitized in the first place. The probability of
a new host of the set $H_{1}-\left\{ 1\right\} $ being parasitized
is then $q-p_{1}$. This value is less than any other value $q-p_{j}$
with $j\neq1$. It is therefore more difficult for a host of the set
$H_{1}-\left\{ 1\right\} $ to be parasitized than for a host of the
set $H_{1}-\left\{ j\right\} $, $j\neq1$, to be parasitized. The
repeated application of this reasoning explains the first inequality
of the proposition. The second inequality can be explained in a similar
manner.

$\:$

For simplicity, we denote $T_{k,N}$ by $T_{k}$ in the particular
case in which $N_{1}=N$. Recalling the definitions of these random
variables and the random variables $Y_{n}$, we obtain the following
relations 
\[
P(Y_{n}\leq k-1)=P(T_{k}>n),
\]
then 
\[
P(Y_{n}\leq k-1)=1-P(T_{k}\leq n)
\]
and 
\[
P(T_{k}=n)=P(Y_{n-1}\leq k-1)-P(Y_{n}\leq k-1).
\]
From (\ref{eq:Pnk notacion conjuntos}), (\ref{eq:FuncDistribYn})
and above equalities we see that 
\begin{equation}
P(T_{k}>n)=\sum_{A\subset H,\left|A\right|\leq k-1}(-1)^{k-1-\left|A\right|}\dbinom{N-\left|A\right|-1}{k-1-\left|A\right|}\: p_{A}^{n},\label{eq:P(Zk>n)}
\end{equation}
\[
P(T_{k}\leq n)=1-\sum_{A\subset H,\left|A\right|\leq k-1}(-1)^{k-1-\left|A\right|}\dbinom{N-\left|A\right|-1}{k-1-\left|A\right|}\: p_{A}^{n}
\]
and 
\[
P(T_{k}=n)=\sum_{A\subset H,\left|A\right|\leq k-1}(-1)^{k-1-\left|A\right|}\dbinom{N-\left|A\right|-1}{k-1-\left|A\right|}\: p_{A}^{n-1}-
\]
\[
\sum_{A\subset H,\left|A\right|\leq k-1}(-1)^{k-1-\left|A\right|}\dbinom{N-\left|A\right|-1}{k-1-\left|A\right|}\: p_{A}^{n}=
\]
\[
\sum_{A\subset H,\left|A\right|\leq k-1}(-1)^{k-1-\left|A\right|}\dbinom{N-\left|A\right|-1}{k-1-\left|A\right|}\: p_{A}^{n-1}(1-p_{A}).
\]

$E(T_{k})=E(T_{k,N})$ is the expected number of draws required for
$k$ hosts are parasitized. Different expressions have been described
for this expectation \citep{BonehHofri,FlajoletBirthdayParadox}.
From (\ref{eq:P(Zk>n)}) it follows immediately that 
\[
E(T_{k})=\sum_{n=0}^{\infty}P(T_{k}>n)=\sum_{n=0}^{\infty}\left(\sum_{A\subset H,\left|A\right|\leq k-1}(-1)^{k-1-\left|A\right|}\dbinom{N-\left|A\right|-1}{k-1-\left|A\right|}\: p_{A}^{n},\right)=
\]
\[
\sum_{A\subset H,\left|A\right|\leq k-1}(-1)^{k-1-\left|A\right|}\dbinom{N-\left|A\right|-1}{k-1-\left|A\right|}\:\frac{1}{1-p_{A}}.
\]
This expression was obtained in \citet{FlajoletBirthdayParadox}.
In \citet{BonehHofri} the following expression was obtained, 
\[
E(T_{k})=\sum_{r=0}^{k-1}\left\Vert u^{r}\right\Vert {\displaystyle \int_{t\geq0}}\prod_{i=1}^{N}(1+u(e^{p_{i}t}-1))e^{-t}\, dt,
\]
where $\left\Vert x^{r}\right\Vert {\displaystyle f(x)}$ is the coefficient
of $x^{r}$ in the power series development of $f(x)$.

$\:$

If $k=N_{1}=N$, then $E(T_{N})=E(T_{N,N})$ is the expected number
of draws required to obtain complete parasitism. From (\ref{eq:E(Zk)}) 

\begin{equation}
E(T_{N})=\sum_{(i_{1}i_{2}...i_{N})\in\Pi_{N}}\left(\sum_{r=0}^{N-1}\dfrac{1}{1-\sum_{j=1}^{r}p_{i_{j}}}\right)\dfrac{\prod_{i=1}^{N}p_{i}}{\prod_{k=1}^{N}\sum_{j=k}^{N}p_{i_{j}}}\label{eq:E(TN in Boneh)-1}
\end{equation}
where $\Pi_{N}$ is the group of permutations of $\left\{ 1,2,...N\right\} .$
This expression for $E(T_{N})$ is proposed in \citet{BonehHofri}.
The authors provide no proof for this formula, and we have found no
proof elsewhere.

\section{Applications to various risk distributions}

In this section, we consider two different risk distributions on the
set of hosts and compute the most relevant values for every each.

$\:$

$\mathbf{The\; uniform\; distribution.}$ The situation in which risk
is distributed uniformly, i.e. all the hosts have the same probability
of being parasitized, with: 
\begin{equation}
p_{1}=p_{2}=...=p_{N}=\dfrac{1}{N}\label{eq:p1=00003D00003Dp2=00003D00003D...=00003D00003DpN}
\end{equation}
has been widely studied. In this case, the expectation and variance
of the random variable $Y_{n}$ representing the number of parasitized
hosts after $n$ draws are 
\[
E(Y_{n})=N-\dfrac{(N-1)^{n}}{N^{n-1}},
\]
\[
Var(Y_{n})=\dfrac{(N-1)(N-2)^{n}}{N^{n-1}}+\dfrac{(N-1)^{n}(N^{n-1}-(N-1)^{n})}{N^{2n-3}}.
\]
and the expected number of draws for $k$ new hosts to be parasitized
(\ref{eq:E(Zk)}) is 
\[
E(T_{k,N_{1}})=N\left(\dfrac{1}{N_{1}}+\dfrac{1}{N_{1}-1}+...+\dfrac{1}{N_{1}-k+1}\right),
\]
which, in the case in which $k=N$, can be written as the following
well-known formula 
\[
E(T_{N})=N\left(1+\dfrac{1}{2}+\dfrac{1}{3}...+\dfrac{1}{N}\right).
\]
It is clear that in this case the upper and lower bounds for $E(T_{k,N_{1}})$
obtained in Proposition 4.3, are both equal to $E(T_{k,N_{1}})$,
and all the probabilities $P(D_{i_{1}i_{2}...i_{k}})$ are equal to
$\dfrac{1}{N_{1}^{(k)}}$.

$\:$

$\mathbf{Two\; kinds\; of\; hosts.}$ The two types of host situation
is an idealization of the following cases. Hosts which are dead, either
because they were previously parasitized or because they produced
artifacts such as mines and galls, remain in the ecosystem for much
longer than the existence of the host. They can make up to 90\% of
the host population. They can be still attractive to parasitoids long
after the host death. Parasitoids will not lay eggs in them, but they
will be checked carefully, implying a waste of time of up to 20\%
\citep{Casas1989,Casasetal2004}. In such cases, it is possible to
envision two categories, living and dead hosts, while being interested
in the rate of parasitism of the living ones only.

Let us now consider the situation in which there are two kinds of
hosts and, therefore, two different probabilities of being detected
by a parasitoid.

In a population of $N$ hosts, each of the hosts 1, 2, ..., $m$ has
a probability $\alpha$ of being parasitized, and each hosts $m+1$,
$m+2$, ..., $N$ has a probability $\beta$ of being parasitized,
such that

\begin{equation}
\begin{array}{c}
p_{1}=p_{2}=...=p_{m}=\alpha,\\
p_{m+1}=p_{m+2}=...=p_{N}=\beta.
\end{array}\label{eq:pi=00003D00003Dalfa pj=00003D00003Dbeta}
\end{equation}

The probability of host 1 being visited $r_{1}$ times, host 2 $r_{2}$
times, etc, for $r_{1}+r_{2}+...+r_{N}=n$, given by (\ref{eq:Prob Sn1=00003D00003Dr1, Sn2=00003D00003Dr2, etc})
is in this case 
\[
P(S_{n1}=r_{1},S_{n2}=r_{2},...,S_{nN}=r_{N})=\frac{n!}{r_{1}!r_{2}!...r_{N}!}{\displaystyle \alpha}^{\sum_{i\leq m}r_{i}}\beta^{\sum_{i>m}r_{i}}
\]

\[
0\leq r_{1}\leq n,0\leq r_{2}\leq n,...,0\leq r_{N}\leq n,\quad r_{1}+r_{2}+...+r_{N}=n.
\]
The probability that, after $n$ draws host $i_{1}$ had been chosen
$r_{i_{1}}$ times by the parasitoids, host $i_{2}$ $r_{i_{2}}$
times and host $i_{h}$ $r_{i_{h}}$ times, without taking the other
hosts into account, is given by (\ref{eq:prob marginal}). It is equal
to 
\[
P(S_{ni_{1}}=r_{i_{1}},S_{ni_{2}}=r_{i_{2}},...,S_{ni_{h}}=r_{i_{h}})=
\]
\[
\frac{n!}{r_{i_{1}}!r_{i_{2}}!...r_{i_{h}}!(n-\sum r_{i_{j}})!}\;\alpha^{\sum_{i_{j}\leq m}r_{i_{j}}}\beta^{\sum_{i_{j}>m}r_{i_{j}}}\left(1-\sum_{i_{j}\leq m}\alpha-\sum_{i_{j}>m}r_{i_{j}}\beta\right)^{n-\sum r_{i_{j}}}.
\]

We will now calculate the expected number of parasitized hosts after
$n$ draws with this risk distribution, using the results obtained
in Section 2.

Let $Y_{n}$ be the random variable representing the number of parasitized
hosts after $n$ draws. From (\ref{eq:Pnk notacion conjuntos}) it
follows that

\[
P(Y_{n}=k)=\sum_{j=1}^{k}(-1)^{k-j}\dbinom{N-j}{k-j}{\displaystyle \sum_{i=0}^{j}}\dbinom{m}{i}\:\dbinom{N-m}{j-i}\:(i\alpha-(j-i)\beta)^{n}
\]
and the expected value of $Y_{n}$, (\ref{eq:E(Yn)}), is equal to
\[
E(Y_{n})=N-m(1-\alpha)^{n}-(N-m)(1-\beta)^{n}.
\]

To compute the expected number of draws for $k$ hosts of a set $H_{1}\subset H$
of unparasitized hosts to be parasitized, we will name the hosts of
the set $H_{1}$, hosts 1, 2, ..., $N_{1}$. Without any loss of generality,
we can assume $p_{1}=p_{2}=...=p_{m_{1}}=\alpha$ and $p_{m_{1}+1}=p_{m_{1}+2}=...=p_{N_{1}}=\beta$.
Let $\Pi_{k}$ be the set of all $k$-permutations of the integers
1, 2, ..., $N_{1}$. For every $I=(i_{1},i_{2},...,i_{k})\in\Pi_{k}$,
let $A_{I}\subset\{1,2,...,k\}$ be the set defined by $j\in A_{I}$
if $i_{j}\leq m_{1}$. It is clear that the probability $P(D_{i_{1},i_{2},...,i_{k}})=P(D_{I})$
given by (\ref{eq:P(Di1i2...in1)}) is, in this case, 
\[
P(D_{I})=\dfrac{\alpha^{\left|A_{I}\right|}\beta^{k-\left|A_{I}\right|}}{p(p-\gamma_{1})(p-\sum_{j=1}^{2}\gamma_{j})...(p-\sum_{j=1}^{k-1}\gamma_{j})},
\]
where 
\begin{equation}
\gamma_{j}=\left\{ \begin{array}{cc}
\alpha & \:\mathrm{if\;}j\in A_{I}\\
\beta & \mathrm{\: if\;}j\notin A_{I}
\end{array}\right.\label{eq:gamma sub j}
\end{equation}
Then, if $A_{I}=A_{I'}$ for $I\in\Pi_{k}$ and $I'\in\Pi_{k}$, it
follows directly that 
\[
P(D_{I})=P(D_{I'}).
\]
We can therefore define an equivalence relation on $\Pi_{k}$ as follows:
$I$ is related to $I'$ if $A_{I}=A_{I'}$. We denote by $\overline{I}$
the equivalence class of $I$, and by $\overline{\Pi}_{k}$ the set
whose elements are the equivalence classes of the elements of $\Pi_{k}$,
that is 
\[
\overline{\Pi}_{k}=\left\{ \overline{I}:I\in\Pi_{k}\right\} .
\]
There are as many equivalence classes as subsets of $\{1,2,...,k\}$
with cardinalities greater than or equal to $\max\left\{ 0,k-n_{1}\right\} $,
where $n_{1}=N_{1}-m_{1}$, and less than or equal to $\min\left\{ k,m_{1}\right\} $,
and the cardinalities of these equivalence classes are 
\[
\left|\bar{I}\right|=m_{1}^{(h)}n_{1}^{(k-h)}\qquad\mathrm{if}\quad\left|A_{I}\right|=h.
\]
Given the above considerations, it is clear that $E(T_{k,N_{1}})$
can be written in this case as: 
\[
\begin{array}{ccc}
E(T_{k,N_{1}}) & = & \sum_{I\in\Pi_{k}}\left(\dfrac{1}{p}+\dfrac{1}{(p-\gamma_{1})}+...+\dfrac{1}{(p-\sum_{j=1}^{k-1}\gamma_{j})}\right)\end{array}=
\]
\[
\dfrac{\alpha^{\left|A_{I}\right|}\beta^{k-\left|A_{I}\right|}}{p(p-\gamma_{1})(p-\sum_{j=1}^{2}\gamma_{j})...(p-\sum_{j=1}^{k-1}\gamma_{j})}=
\]
\[
\sum_{\bar{I}\in\bar{\Pi}_{k}}\sum_{I\in\bar{I}}\left(\dfrac{1}{p}+\dfrac{1}{(p-\gamma_{1})}+...+\dfrac{1}{(p-\sum_{j=1}^{k-1}\gamma_{j})}\right)
\]
\[
\dfrac{\alpha^{\left|A_{I}\right|}\beta^{k-\left|A_{I}\right|}}{p(p-\gamma_{1})(p-\sum_{j=1}^{2}\gamma_{j})...(p-\sum_{j=1}^{k-1}\gamma_{j})}=
\]
\[
\sum_{\bar{I}\in\bar{\Pi}_{k}}m_{1}^{(\left|A_{I}\right|)}n_{1}^{(k-\left|A_{I}\right|)}\left(\dfrac{1}{p}+\dfrac{1}{(p-\gamma_{1})}+...+\dfrac{1}{(p-\sum_{j=1}^{k-1}\gamma_{j})}\right)
\]
\[
\dfrac{\alpha^{\left|A_{I}\right|}\beta^{k-\left|A_{I}\right|}}{p(p-\gamma_{1})(p-\sum_{j=1}^{2}\gamma_{j})...(p-\sum_{j=1}^{k-1}\gamma_{j})}=
\]
\[
{\displaystyle \sum_{h=\max\left\{ 0,k-n_{1}\right\} }^{\min\left\{ k,m_{1}\right\} }}\sum_{\left|A_{I}\right|=h}m_{1}^{(h)}n_{1}^{(k-h)}\left(\dfrac{1}{p}+\dfrac{1}{(p-\gamma_{1})}+...+\dfrac{1}{(p-\sum_{j=1}^{k-1}\gamma_{j})}\right)
\]
\[
\dfrac{\alpha^{h}\beta^{k-h}}{p(p-\gamma_{1})(p-\sum_{j=1}^{2}\gamma_{j})...(p-\sum_{j=1}^{k-1}\gamma_{j})}.
\]
where $\gamma_{j}$ is defined by (\ref{eq:gamma sub j}).

Let us suppose that 
\[
\alpha>\beta.
\]

To obtain an upper bound for $E(T_{k,N_{1}})$, we distinguish two
cases, $k\leq m_{1}$ and $k>m_{1}$. If $k\leq m_{1}$ then

\[
E(T_{k,N_{1}}|D_{1,2,...,k})=
\]

\[
\dfrac{1}{m_{1}\alpha+n_{1}\beta}+\dfrac{1}{(m_{1}-1)\alpha+n_{1}\beta}+...+\dfrac{1}{(m_{1}-k+1)\alpha+n_{1}\beta},
\]
if $k>m_{1}$, this upper bound is

\[
E(T_{k,N_{1}}|D_{1,2,...,k})=
\]
\[
\dfrac{1}{m_{1}\alpha+n_{1}\beta}+\dfrac{1}{(m_{1}-1)\alpha+n_{1}\beta}+...+\dfrac{1}{n_{1}\beta}+
\]
\[
\dfrac{1}{(n_{1}-1)\beta}+...+\dfrac{1}{(n_{1}+m_{1}-k+1)\beta}.
\]
Similarly, to obtain a lower bound for $E(T_{k,N_{1}})$ we distinguish
the cases $k\leq n_{1}$ and $k>n_{1}$. If $k\leq n_{1}$ this lower
bound is 
\[
E(T_{k,N_{1}}|D_{N_{1},N_{1}-1,...,N_{1}-k+1})=
\]

\[
\dfrac{1}{m_{1}\alpha+n_{1}\beta}+\dfrac{1}{m_{1}\alpha+(n_{1}-1)\beta}+...+\dfrac{1}{m_{1}\alpha+(n_{1}-k+1)\beta},
\]
and if $k>n_{1}$, a lower bound for $E(T_{k,N_{1}})$ is 
\[
E(T_{k,N_{1}}|D_{N_{1},N_{1}-1,...,N_{1}-k+1})=
\]

\[
\dfrac{1}{m_{1}\alpha+n_{1}\beta}+\dfrac{1}{m_{1}\alpha+(n_{1}-1)\beta}+...+\dfrac{1}{m_{1}\alpha}+
\]
\[
\dfrac{1}{(m_{1}-1)\alpha}+\dfrac{1}{(n_{1}+m_{1}-k+1)\alpha}.
\]
The maximum of the values $P(D_{i_{1},i_{2},...,i_{k}})$ is 
\[
P(D_{1,2,...,k})=\left\{ \begin{array}{ccc}
\dfrac{\alpha^{k}}{\prod_{h=0}^{k-1}((m_{1}-h)\alpha+n_{1}\beta)}, & \mathrm{if} & k\leq m_{1}\\
\dfrac{\alpha^{m_{1}}\beta^{k-m_{1}}}{\prod_{h=0}^{m_{1}}((m_{1}-h)\alpha+n_{1}\beta)\prod_{l=1}^{k-m_{1}-1}(n_{1}-l)\beta}, & \mathrm{if} & k>m_{1}
\end{array}\right.
\]

$\:$

In the extreme case that there is only one host with a probability
$\alpha$ of being parasitized and the others have probability $\beta$
of being parasitized, we obtain the following expressions for $E(T_{k,N_{1}})$.

If the host with probability $\alpha$ of being parasitized does not
belong to set $H_{1}$, then 
\[
E(T_{k,N_{1}})=\dfrac{1}{\beta}\sum_{j=0}^{k-1}\dfrac{1}{N_{1}-h}.
\]
If the host with probability $\alpha$ of being parasitized belongs
to set $H_{1}$, then

\[
E(T_{k,N_{1}})=(N_{1}-1)^{(k)}\left(\dfrac{1}{\alpha+(N_{1}-1)\beta}+\dfrac{1}{\alpha+(N_{1}-2)\beta}+...+\dfrac{1}{\alpha+(N_{1}-k)\beta}\right)
\]

\[
\dfrac{\beta^{k}}{(\alpha+(N_{1}-1)\beta)(\alpha+(N_{1}-2)\beta)...(\alpha+(N_{1}-k)\beta)}+
\]

\[
(N_{1}-1)^{(k-1)}\sum_{j=1}^{k-1}\left(\dfrac{1}{\alpha+(N_{1}-1)\beta}+\dfrac{1}{\alpha+(N_{1}-2)\beta}+...+\dfrac{1}{\alpha+(N_{1}-j)\beta}+\right.
\]
\[
\left.\dfrac{1}{(N_{1}-j)\beta}+\dfrac{1}{(N_{1}-j-1)\beta}+...+\dfrac{1}{(N_{1}-k+1)\beta}\right)
\]

\[
\dfrac{\alpha\beta^{k-1}}{(\alpha+(N_{1}-1)\beta)(\alpha+(N_{1}-2)\beta)...(\alpha+(N_{1}-j)\beta)(N_{1}-j)\beta(N_{1}-j-1)\beta...(N_{1}-k+1)\beta}+
\]
\[
(N_{1}-1)^{(k-1)}\left(\dfrac{1}{\alpha+(N_{1}-1)\beta}+\dfrac{1}{\alpha+(N_{1}-2)\beta}+...+\dfrac{1}{\alpha+(N_{1}-k)\beta}\right)
\]
\[
\dfrac{\alpha\beta^{k-1}}{(\alpha+(N_{1}-1)\beta)(\alpha+(N_{1}-2)\beta)...(\alpha+(N_{1}-k)\beta)}.
\]

$\:$

\section{Relationship between the risk distribution and the speed of parasitization}

In the preceding sections, we studied the process of parasitization
for a given risk distribution in the set of hosts. In this section
we compare this process for different risk distributions. We show
how parasitization speed depends on the risk distribution, and its
scatter in particular. We use the concept of ``majorization'' to
formalize the idea that risk distributions have different degrees
of spread. This notion dates from the start of the 20th century. A
comprehensive review of the theory can be found in \citet{Majorization}.

Less spread distributions are associated with faster parasitization.
In other words, the more spread out the risk distribution, the larger
the number of draws required for a given number of hosts to be parasitized.
Thus the distribution function for the first time parasitization of
a given number of hosts, viewed as a function of the vector $p$,
is Schur convex (see the definition at the end of this section). The
mathematical community studying the coupon collector problem seems
to be largely unaware of it, but this result is not new and can be
found in \citet{WongYue1973WeakDominance}. This result constitutes
the first part of Theorem 6.1. We give a proof more concise and clearer
than previous proposal. Moreover, our method provides a precise result
for strict Schur convexity. This refinement constitutes the second
part of Theorem 6.1. We make use in our proof of the relationship
between the concept of majorization and the numerical operation known
as \textquotedbl{}Robin Hood transfer\textquotedbl{}, described below.

$\:$

In this section, we work with different risk distributions, requiring
further notation and definitions. Given a risk distribution $p=(p_{1},p_{2},...,p_{N})$,
we denote by $P_{p}$ the probability distribution induced by $p$
on the $\sigma$-field over the space of the all the possible incidences
of parasitization.

Given $(p_{1},p_{2},...,p_{N})$ in $\mathbb{R}^{N}$ , we denote
by $(p_{\bar{1}},p_{\bar{2}},...,p_{\bar{N}})$ the $N$-uple obtained
by permutation of $p_{i}$ such that $p_{\bar{1}}\geq p_{\bar{2}}\geq\ldots\geq p_{\bar{N}}$.

The following definitions are given in \citet{Majorization}.

$\:$

$\mathbf{Definition\:6.1.}$ Let $p_{1},p_{2},...,p_{N}$, $q_{1},q{}_{2},...,q{}_{N}$,
be real numbers. We say that $p=(p_{1},p_{2},...,p_{N})$ is majorized
by $q=(q_{1},q{}_{2},...,q{}_{N})$, and we write $p\prec q$, if
\[
\sum_{i=1}^{k}p_{\bar{i}}\leq\sum_{i=1}^{k}q{}_{\bar{i}}\quad\mathrm{for}\quad i=1,2,...,N-1
\]
and 
\[
\sum_{i=1}^{N}p_{\bar{i}}=\sum_{i=1}^{N}q{}_{\bar{i}}.
\]

$\:$

It is clear that when we apply this definition to the comparison of
two risk distributions, the last equality is trivially satisfied.

$\:$

Let $q=(q_{1},q{}_{2},...,q{}_{N})\in\mathbb{R}^{N}$. If $q_{h}<q_{k}$
we can transfer an amount $\Delta$, $0<\Delta<q_{k}-q_{h}$ from
$q_{k}$ to $q_{h}$ to obtain the following new risk distribution
$q'=(q'_{1}$, $q'_{2}$, ...,$q'_{N})$, where $q'_{h}=q_{h}+\varDelta$,
$q'_{k}=q_{k}-\varDelta$ and $q'_{i}=q_{i}$ for $i\neq h,k$. Then,
$q'$ is less spread out than the initial distribution, that is, $q'\prec q$.
Such operations involving the shifting of some ``income'' from one
individual to a poorer individual, are described, somewhat poetically,
as Robin Hood transfers \citep{ArnoldRobinHoodTransfer}. If we define
$\alpha=1-\tfrac{\Delta}{q_{k}-q_{h}}$ then we can write $q'_{h}=q_{h}+\varDelta=\alpha q_{h}+(1-\alpha)q_{k}$
and $q'_{k}=q_{k}-\varDelta=\alpha q_{k}+(1-\alpha)q_{h}$.

$\:$

$\mathbf{Proposition\:6.1.}$ The following conditions are equivalent:

a) $p\prec q$, 

b) $p$ can be derived from $q$ by successive applications of a finite
number of Robin Hood transfers.

$\:$

It is not difficult to prove this equivalence. It was proved for the
first time in \citet{Muirhead1903} for vectors of non-negative integer
components.

$\:$

$\mathbf{Lemma\:6.1.}$ Let $k$ and $N$ be integers satisfying $1<k\leq N-1$,
then 
\[
\sum_{0\leq r\leq k-1}(-1)^{k-1-r}\dbinom{N-2}{r}\:\dbinom{N-2-r}{k-1-r}=0.
\]

$\mathbf{Proof.}$ For $x\in\mathbb{R}$ the equality 
\[
\sum_{0\leq r\leq n}\dbinom{a}{r}\;\dbinom{x}{n-r}=\dbinom{a+x}{n},
\]
is satisfied, where $\dbinom{x}{h}$ is defined by (\ref{eq:NumeroCombinatorio})
then 
\[
\sum_{0\leq r\leq n}\dbinom{a}{r}\:\dbinom{-b}{n-r}=\dbinom{a-b}{n},
\]
and 
\[
\sum_{0\leq r\leq n}(-1)^{n-r}\dbinom{a}{r}\:\dbinom{b+n-r-1}{n-r}=\sum_{0\leq r\leq n}\dbinom{a}{r}\:\dbinom{-b}{n-r}=\dbinom{a-b}{n}.
\]
Then, we obtain, with $a=N-2$, $n=k-1$, $b=N-k$ 
\[
\sum_{0\leq r\leq k-1}(-1)^{k-1-r}\dbinom{N-2}{r}\:\dbinom{N-2-r}{k-1-r}=\dbinom{k-2}{k-1}=0,
\]
and the lemma follows.

$\:$

$\mathbf{Lemma\:6.2}.$ Let $q_{1}$, $q_{2}$, ...,$q_{M}$ be non-negative
real numbers and $I=\left\{ 1,2,...,M\right\} $. For every $A\subset I$
let $q_{A}=\sum_{i\in A}q_{i}$. Then, for any integer $m\geq0$ ,
\[
\sum_{A\subset I,\left|A\right|\leq r}(-1)^{r-\left|A\right|}\dbinom{M-\left|A\right|}{r-\left|A\right|}\: q_{A}^{m}\geq0.
\]

Moreover, if $m\geq r$ and at least $r$ of the values $q_{1}$,
$q_{2}$, ...,$q_{M}$ are greater than zero, then 
\[
\sum_{A\subset I,\left|A\right|\leq r}(-1)^{r-\left|A\right|}\dbinom{M-\left|A\right|}{r-\left|A\right|}\: q_{A}^{m}>0.
\]

$\:$

$\mathbf{Proof.}$ If all the $q_{i}$ are zero, there is nothing
to prove. Let us suppose that $s=\sum_{i=1}^{M}q_{i}>0$. Let $p_{i}=q_{i}/s$,
$i=1,2,...,M$. These values define the probability distribution $p=(p_{1},p_{2},...,p_{M})$
on $I$. From (\ref{eq:Pnk notacion conjuntos}) it follows that 
\[
\sum_{A\subset I,\left|A\right|\leq r}(-1)^{r-\left|A\right|}\dbinom{M-\left|A\right|}{r-\left|A\right|}\: q_{A}^{m}=
\]
\[
s^{m}\sum_{A\subset I,\left|A\right|\leq r}(-1)^{r-\left|A\right|}\dbinom{M-\left|A\right|}{r-\left|A\right|}\: p_{A}^{m}=s^{m}P_{p}(Y_{m}=r),
\]
which proves the lemma.

$\:$

Let $p=(p_{1},p_{2},...,p_{N})$ denote a probability distribution
$p$ over the set $H$. Suppose that $p$ is not uniform. We can assume
$p_{1}<p_{2}$ without loss of generality. Let $0<h\leq\tfrac{p_{2}-p_{1}}{2}$,
$\alpha=1-\tfrac{h}{p_{2}-p_{1}}$. We then define a new risk distribution
$p'$ by applying a Robin Hood transfer as follows 
\begin{equation}
p'=(p_{1}+h,p_{2}-h,p_{3},p_{4},...,p_{N})=(\alpha p_{1}+(1-\alpha)p_{2},\alpha p_{2}+(1-\alpha)p_{1},p_{3},...,p_{N}).\label{eq:p' defined from p}
\end{equation}
 We indeed have $p'\prec p$.

\textcolor{red}{$\:$}

$\mathbf{Theorem\:6.1.}$ Let $p$ be a non uniform probability distribution
over $H$. Without loss of generality, we can assume that $p_{1}<p_{2}$.
Let $p'$ be defined by (\ref{eq:p' defined from p}). Then, for all
$k$ between 1 and $N-1$, 
\begin{equation}
P_{p}(Y_{n}\leq k)\geq P_{p'}(Y_{n}\leq k),\label{eq:Pp>=00003D00003DPp'}
\end{equation}
which is equivalent to

\begin{equation}
P_{p}(T_{k+1}\leq n)\leq P_{p'}(T_{k+1}\leq n)\label{eq:Tk+1p<=00003D00003DTk+1p'}
\end{equation}
Moreover, if at least $k-1$ of the values $p_{3},p_{4},...,p_{N}$
are non-zero, then

\begin{equation}
P_{p}(Y_{n}\leq k)>P_{p'}(Y_{n}\leq k),\quad n=k+1,k+2,k+3...\label{eq:Pp>Pp'}
\end{equation}
which is equivalent to 
\begin{equation}
P_{p}(T_{k+1}\leq n)<P_{p'}(T_{k+1}\leq n),\quad n=k+1,k+2,k+3...\label{eq:Tk+1p<Tk+1p'}
\end{equation}
where $p'$ is defined by (\ref{eq:p' defined from p}).

$\:$

$\mathbf{Proof.}$ Let $H'=\left\{ 3,4,...,N\right\} $. According
to (\ref{eq:FuncDistribYn}) we have: 
\[
P_{p}(Y_{n}\leq k)=\sum_{A\subset H,\left|A\right|\leq k}(-1)^{k-\left|A\right|}\dbinom{N-\left|A\right|-1}{k-\left|A\right|}\: p_{A}^{n}=
\]
\[
\sum_{A\subset H',\left|A\right|\leq k}(-1)^{k-\left|A\right|}\dbinom{N-\left|A\right|-1}{k-\left|A\right|}\: p_{A}^{n}+
\]
\[
\sum_{A\subset H',\left|A\right|\leq k-1}(-1)^{k-1-\left|A\right|}\dbinom{N-2-\left|A\right|}{k-1-\left|A\right|}\:(p_{A}+p_{1})^{n}+
\]
\[
\sum_{A\subset H',\left|A\right|\leq k-1}(-1)^{k-1-\left|A\right|}\dbinom{N-2-\left|A\right|}{k-1-\left|A\right|}\:(p_{A}+p_{2})^{n}+
\]
\[
\sum_{A\subset H',\left|A\right|\leq k-2}(-1)^{k-2-\left|A\right|)}\dbinom{N-3-\left|A\right|}{k-2-\left|A\right|}\:(p_{A}+p_{1}+p_{2})^{n}.
\]

Then 
\[
P_{p}(Y_{n}\leq k)-P_{p'}(Y_{n}\leq k)=
\]
\[
\sum_{A\subset H',\left|A\right|\leq k-1}(-1)^{k-1-\left|A\right|}\dbinom{N-2-\left|A\right|}{k-1-\left|A\right|}\:(p_{A}+p_{1})^{n}+
\]
\[
\sum_{A\subset H',\left|A\right|\leq k-1}(-1)^{k-1-\left|A\right|}\dbinom{N-2-\left|A\right|}{k-1-\left|A\right|}\:(p_{A}+p_{2})^{n}-
\]
\[
\sum_{A\subset H',\left|A\right|\leq k-1}(-1)^{k-1-\left|A\right|}\dbinom{N-2-\left|A\right|}{k-1-\left|A\right|}\:(p_{A}+p_{1}+h)^{n}-
\]
\[
\sum_{A\subset H',\left|A\right|\leq k-1}(-1)^{k-1-\left|A\right|}\dbinom{N-2-\left|A\right|}{k-1-\left|A\right|}\:(p_{A}+p_{2}-h)^{n}.
\]
Let $f$ be the real function defined by 
\[
f(x)=\sum_{A\subset H',\left|A\right|\leq k-1}(-1)^{k-1-\left|A\right|}\:\dbinom{N-2-\left|A\right|}{k-1-\left|A\right|}\:(p_{A}+x)^{n},\qquad x\in\mathbf{\mathbb{\mathbf{R}}}
\]
This function is a polynomial of degree less than or equal to $n$.
The coefficient of $x^{n}$ is equal to 
\[
\sum_{A\subset H',\left|A\right|\leq k-1}(-1)^{k-1-\left|A\right|}\dbinom{N-2-\left|A\right|}{k-1-\left|A\right|}=
\]
\[
\sum_{0\leq r\leq k-1}(-1)^{k-1-r}\dbinom{N-2}{r}\:\dbinom{N-2-r}{k-1-r}
\]
and this is equal to 0 by Lemma 6.1. The coefficient of $x^{n-j}$
for $j=1,2,\:...,n$ is 
\[
\dbinom{n}{j}\sum_{A\subset H',\left|A\right|\leq k-1}(-1)^{k-1-\left|A\right|}\dbinom{N-2-\left|A\right|}{k-1-\left|A\right|}\: p_{A}^{j},
\]
by the first part of Lemma 6.2 with $I=H'=\left\{ 3,4,...,N\right\} $,
$M=\left|I\right|=N-2$, $m=j$ and $r=k-1$, it follows that these
coefficients are greater than or equal to zero. This polynomial function
is then convex on $\left[0,+\infty\right)$, so that 
\[
f(p_{1})+f(p_{2})\geq f(\alpha p_{1}+(1-\alpha)p_{2})+f(\alpha p_{2}+(1-\alpha)p_{1})=f(p_{1}+h)+f(p_{2}-h),
\]
\[
f(p_{1})+f(p_{2})-f(p_{1}+h)-f(p_{2}-h)\geq0.
\]
However, this inequality is the same as 
\[
P_{p}(Y_{n}\leq k)-P_{p'}(Y_{n}\leq k)\geq0,
\]
which gives (\ref{eq:Pp>=00003D00003DPp'}). Recalling the relationship
between the random variables $Y_{i}$ and the random variables $T_{j}$,
we also obtain 
\[
P_{p}(T_{k+1}\leq n)\leq P_{p'}(T_{k+1}\leq n),
\]
which is (\ref{eq:Tk+1p<=00003D00003DTk+1p'}).

Moreover, from the second part of Lemma 6.2. it follows that if at
least $k-1$ of the values $p_{3},p_{4},...,p_{N}$ are greater than
zero and $n\geq k+1$, then the coefficient of $x^{n-k+1}$ is greater
than zero, where $n-k+1\geq2$. So, at least one monomial of degree
greater than or equal to 2 appears in the polynomial. The convexity
is then strict, and we can write 
\[
f(p_{1})+f(p_{2})-f(p_{1}+h)-f(p_{2}-h)>0,
\]
and 
\[
P_{p}(Y_{n}\leq k)-P_{p'}(Y_{n}\leq k)>0,\quad n=k+1,k+2,k+3...
\]
which is equivalent to 
\[
P_{p}(Y_{n}\leq k)>P_{p'}(Y_{n}\leq k),\quad n=k+1,k+2,k+3...
\]
and therefore to 
\[
P_{p}(T_{k+1}\leq n)<P_{p'}(T_{k+1}\leq n),\quad n=k+1,k+2,k+3...
\]
This completes the proof.

$\:$

We can state the following corollaries.

$\:$

$\mathbf{Corollary\:6.1}.$ Let $p=(p_{1},p_{2},...,p_{N})$ and $q=(q_{1},q{}_{2},...,q{}_{N})$
be risk distributions on $H=\left\{ 1,2,...,N\right\} $. If $p\prec q$
then, for every $n\geq1$ and every $k\geq1$ 

\begin{equation}
P_{p}(Y_{n}\leq k)\leq P_{q}(Y_{n}\leq k),\label{eq:Pp<=00003D00003DPq}
\end{equation}
is satisfied and

\[
P_{p}(T_{k+1}\leq n)\geq P_{q}(T_{k+1}\leq n).
\]

Furthermore, if the distributions $p$ and $q$ are actually different,
meaning that they do not differ only by a permutation, then the preceding
inequalities are strict, except in trivial cases. More precisely,
denoting by $j$ the number of non zero $p_{i}$ values (and remarking
that the number of non-zero $q_{i}$ values is at most $j$), we have: 
\begin{itemize}
\item If $k\geq n$ or $k\geq j$ then 
\[
P_{p}(Y_{n}\leq k)=P_{q}(Y_{n}\leq k)=1\text{ and }P_{p}(T_{k+1}\leq n)=P_{q}(T_{k+1}\leq n)=0;
\]

\item If $n\geq2$, $k<n$ and $k<j$, then 
\[
P_{p}(Y_{n}\leq k)<P_{q}(Y_{n}\leq k)\text{ and }P_{p}(T_{k+1}\leq n)>P_{q}(T_{k+1}\leq n).
\]

\end{itemize}
$\mathbf{Proof.}$ As it is possible to go from vector $q$ to vector
$p$ by a finite sequence of Robin Hood transfers, the corollary follows
directly from Theorem 6.1, which proves that each transfer decreases
the quantity $P_{p}(Y_{n}\leq k)$. We just have to consider the cases
in which this quantity is strictly decreased.

$\mathbf{Remark\:6.1.}$ We can interpret the results obtained above
in terms of the theory of Schur-convex functions. A real-valued function
$\phi$ defined on a set $\mathcal{A\subset\mathbf{R}}^{\mathit{N}}$
is said to be Schur-convex on $\mathcal{A}$ if, for every $x$ and
$y$ pair of elements in $\mathcal{A}$ such that $x\prec y$ the
inequality $\phi(x)\leq\phi(y)$ holds. The first part of Corollary
6.1 states that the map $p\rightarrow P_{p}(Y_{n}\leq k)$ is Schur-convex.
This was already proved in \citet{WongYue1973WeakDominance}, and
was stated as a conjecture in \citet{Anceaume2015}.

$\:$

$\mathbf{Corollary\:6.2.}$ Let $u=(1/N,1/N,...,1/N)$ be the uniform
distribution on $H=\left\{ 1,2,...,N\right\} $ and $p=(p_{1},p_{2},...,p_{N})$
any other risk distribution on $H$. Then 
\[
P_{u}(Y_{n}\leq k)<P_{q}(Y_{n}\leq k),\quad k=1,2,...,N-1,\quad n=k+1,k+2,...
\]
\[
P_{u}(T_{k+1}\leq n)>P_{q}(T_{k+1}\leq n),\quad k=1,2,...,N-1,\quad n=k+1,k+2,...
\]

$\mathbf{Proof.}$ It can be clearly seen that $u=(1/N,1/N,...,1/N)$
is majorized by any other distribution on $H$ and the corollary follows.

$\:$

$\mathbf{Remark\:6.2.}$ The results obtained in Corollary 6.1 and
Corollary 6.2 can be expressed in terms of a comparison of probability
distributions as follows. If $p\prec q$, then relation (\ref{eq:Pp<=00003D00003DPq})
proves that the random variable $Y_{n}$ defined on the probability
space determined by $p$ on the space of the random sets of $H=\left\{ 1,2,...,N\right\} $
is weakly stochastically dominated by the random variable $Y_{n}$
defined on the probability space determined by $q$. Corollary 6.2
proves that the random variable $Y_{n}$ defined on the probability
space determined by the uniform distribution $u=(1/N,1/N,...,1/N)$
is always weakly stochastically dominated by the random variable $Y_{n}$
defined on the probability space determined by any other probability
distribution on $H$.

$\:$

$\mathbf{Remark\:6.3.}$ After the redaction of this section, we have
seen a similar study in \citet{Anceaume2016}. In particular, they
prove inequalities (30) and (31) of Theorem 6.1. However, our contribution
still presents a real interest, thanks to the quality of the argument
based on use of fundamental formulas (6) and (7) in different contexts,
and because we obtain cases of strict inequalities.

\section{Illustrative examples}

In this section we show graphically the relationships satisfied among
the distribution functions of random variables $Y_{n}$ as well as
the distribution functions of random variables $T_{k}$, when their
corresponding risk distributions are able to be compared by majorization.

The distribution functions of five variables $Y_{n}$ are represented
in graphic A of Figure 2.
\begin{figure}

\includegraphics[width=12cm]{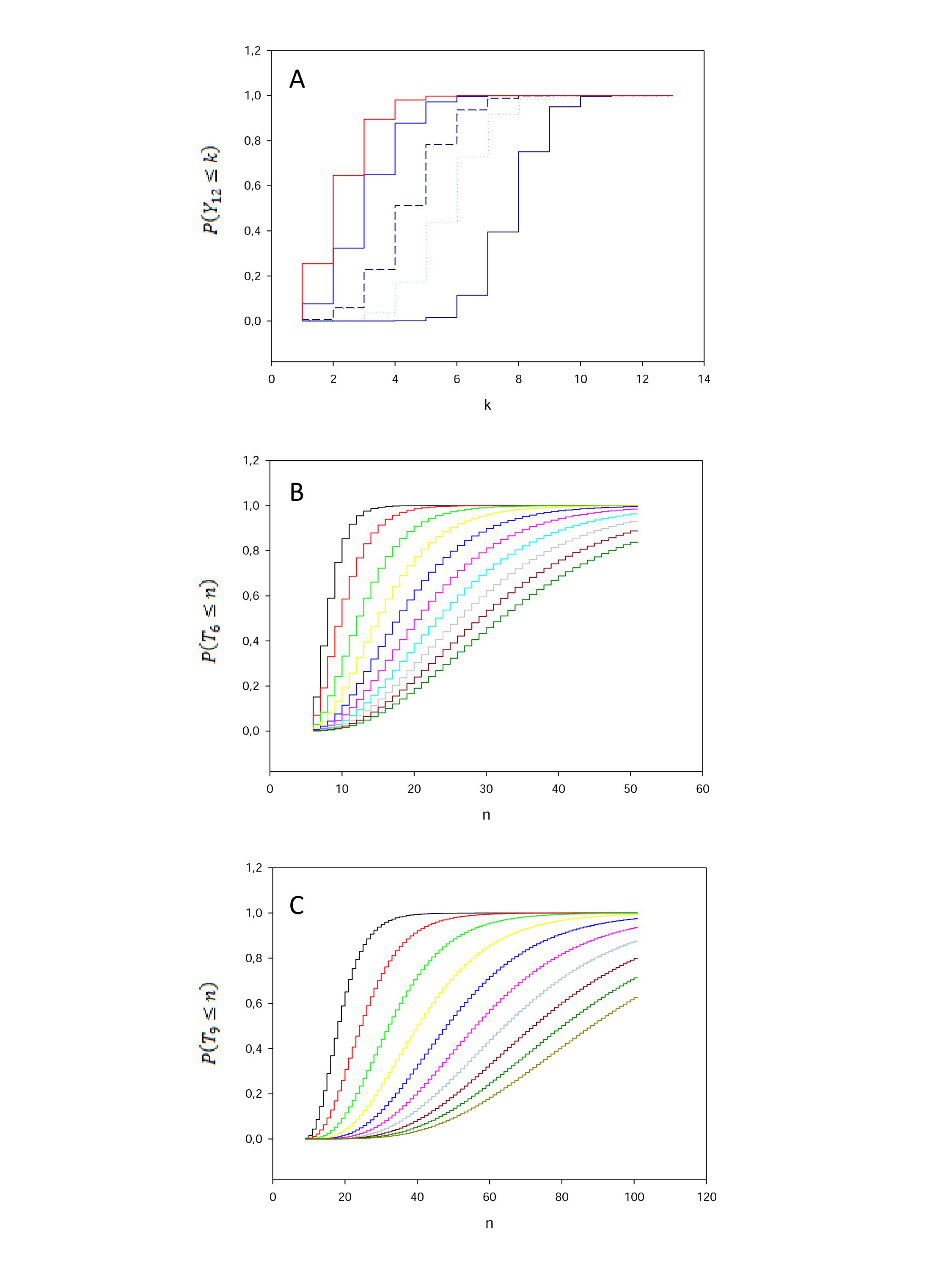}

Figure 2. Graphic A: Distribution functions of five variables $Y_{12}$
corresponding to five different risk distributions, $p_{1},$ $p_{2},$...,
$p_{5},$ satisfying $p_{1}\prec$ $p_{2}\prec$ $p_{3}\prec$ $p_{4}\prec$$p_{5}$.
It can be observed that $P_{p_{i}}(Y_{12}\leq k)<P_{p_{i+1}}(Y_{12}\leq k),\mathrm{for}\quad k=1,2,...,11,\: i=1,2,3,4$.
Graphic B: Distribution functions of ten variables $T_{6}$ corresponding
to ten risk distributions, $p_{1},$ $p_{2},$ ..., $p_{10}$ satisfying
$p_{1}\prec$ $p_{2}\prec$ ... $p_{9}\prec$$p_{10}$. Graphic C:
Distribution functions of ten variables $T_{9}$ corresponding to
the same previous risk distributions. In graphics B and C it can be
observed that $P_{p_{i}}(T_{k}\leq n)>P_{p_{i+1}}(T_{k}\leq n)$,
for $n=k,k+1,...\: i=1,2,...9$
\end{figure}
 They correspond to five different risk distributions, $p_{1},$ $p_{2},$
$p_{3},$ $p_{4},$ and $p_{5},$ satisfying $p_{1}\prec$ $p_{2}\prec$
$p_{3}\prec$ $p_{4}\prec$$p_{5}$. These are distributions on the
set $\left\{ 1,2,...,12\right\} $ (so $N=12$), $p_{1}$ is the uniform
distribution, $p_{i}=(\nicefrac{1}{10i+2},...,\nicefrac{1}{10i+2},\nicefrac{10(i-1)+1}{10i+2})$
for $i=2$ and 3, and $p_{i}=(\nicefrac{1}{45(i-3)+12},...,$ $\nicefrac{1}{45(i-3)+12},\nicefrac{1}{45(i-3)+12},\nicefrac{45(i-3)+1}{45(i-3)+12})$
for $i=4$ and 5. We have also used $n=12$, and it can be observed
that $P_{p_{i}}(Y_{12}\leq k)<P_{p_{i+1}}(Y_{12}\leq k),\mathrm{\; with}\quad k=1,2,...,11,\: i=1,2,3,4$. 

The distribution functions of ten variables $T_{k}$ are represented
in every one of the graphics B and C in Figure 2. $N=10$ and the
risk distributions are the same in both cases; $p_{1}$ is the uniform
distribution and $p_{i}=(\nicefrac{1}{5(i+1)},...,\nicefrac{1}{5(i+1)},$
$\nicefrac{i}{5(i+1)},\nicefrac{4(i-1)+1}{5(i+1)})$, for $i=2,3,...,10$.
For these risk distributions $p_{1}\prec$ $p_{2}\prec$ ... $p_{9}\prec$$p_{10}$
is satisfied. In graphic B of Figure 2, $k=6$ and the values of $n$
lie between 6 to 50. In graphic C of Figure 2, $k=9$ and the values
of $n$ lie from 9 to 100. It can be seen that $P_{p_{i}}(T_{k}\leq n)>P_{p_{i+1}}(T_{k}\leq n)$,
for $n=k,k+1,...\: i=1,2,...9$, in both graphics.

Figure 3 compares distribution functions of random variables $T_{k}$
corresponding to two unrelated risk distributions $p$ and $q$, i.e.
neither $p\prec q$ nor $q\prec p$. 
\begin{figure}

\includegraphics[width=12cm]{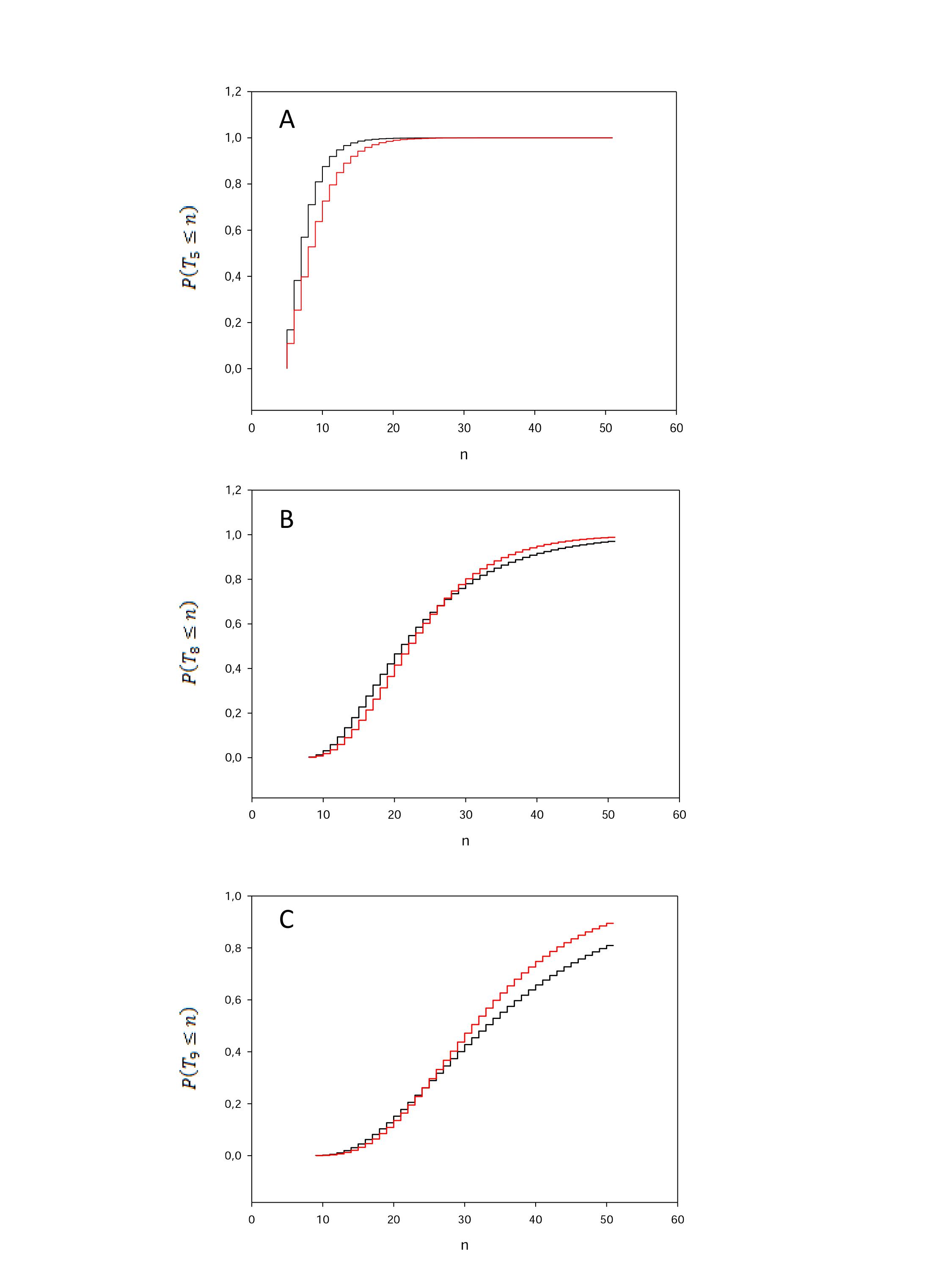}

Figure 3. Comparison of distribution functions of random variables
$T_{k}$ corresponding to two unrelated risk distributions, i.e. neither
$p\prec q$ nor $q\prec p$, to show how these distribution functions
act in different ways depending on the value of $k$. 
\end{figure}
Thus, these distribution functions act in different ways depending
on the value of $k$. We include three different graphics, each bearing
two curves. These curves are the distribution functions of two random
variables $T_{k}$. The risk distributions associated with these random
variables are, in the three graphics, $p=(\nicefrac{3}{85},\nicefrac{3}{85},\nicefrac{3}{85},\nicefrac{3}{85},\nicefrac{3}{85},\nicefrac{12}{85},$
$\nicefrac{13}{85},\nicefrac{13}{85},\nicefrac{13}{85},\nicefrac{19}{85})$
and $q=(\nicefrac{3}{81},\nicefrac{4}{81},\nicefrac{4}{81},\nicefrac{4}{81},$
$\nicefrac{4}{81},\nicefrac{5}{81},\nicefrac{5}{81},\nicefrac{5}{81},\nicefrac{15}{81},\nicefrac{32}{81})$.
In the first graphic $k=5$, in the second $k=8$ and in the third
$k=9$. In the last two cases the distribution functions cross. They
do not cross in the first.

\section{A conjecture on strong dominance}

In Section 6 we used an order relationship between random variables
(or more precisely between their distributions) that can be defined
formally as follows.

$\:$

$\mathbf{Definition\:8.1.}$\emph{ }\textit{\emph{Let $X$ and $X\lyxmathsym{\textasciiacute}$
be two real random variables, defined on probability spaces $\left(\Omega,P\right)$
and $\left(\Omega',P'\right)$, respectively. We say that the random
variable $X$ weakly stochastically dominates the random variable
$X\lyxmathsym{\textasciiacute}$ if the cumulative distribution function
of $X\lyxmathsym{\textasciiacute}$dominates the cumulative distribution
function of $X$, that is, for any }}\textit{$t\in\mathbb{R}$, }

\[
P\left(X\leq t\right)\leq P\text{\textasciiacute}\left(X\text{\textasciiacute}\leq t\right).
\]

The main result of Section 6 is that if $p\prec q$, then the random
variable $Y_{n}$ defined on the probability space $\left(\Omega,P_{p}\right)$
weakly stochastically dominates the random variable $Y_{n}$ defined
on the probability space $\left(\Omega,P_{q}\right).$

A particular case of weak dominance is that one in which inequalities
apply not only to the cumulative distribution functions, but also
to the distributions themselves. We will refer to this situation as
strong dominance, and we provide a formal definition of strong dominance
below, for the case of discrete random variables. (A similar definition
can be given for continuous random variables with densities). In short,
$X$ strongly dominates $X\text{\textasciiacute}$ if, for any small
enough value \textit{d, $P\left(X=d\right)\leq P\text{\textasciiacute}\left(X\text{\textasciiacute}=d\right)$,
}and if for any other possible value\textit{ e, $P\left(X=e\right)\geq P\text{\textasciiacute}\left(X\text{\textasciiacute}=e\right).$} 

$\:$

$\mathbf{Definition\:8.2.}$\emph{ }\textit{\emph{Let $X$ and $X\lyxmathsym{\textasciiacute}$
be two real random variables, defined on probabilities spaces $\left(\Omega,P\right)$
and $\left(\Omega',P'\right)$, respectively, and taking values in
a denumerable set $D$. We say that the random variable $X$ strongly
stochastically dominates the random variable $X\lyxmathsym{\textasciiacute}$
if there is a critical value $c\in\mathbb{R}$ such that, for any}}\textit{
$d\in D$}

\textit{\emph{$\bullet$if }}\emph{$d\leq c$,}\textit{\emph{ then}}\emph{
$P\left(X=d\right)\leq P\text{\textasciiacute}\left(X\text{\textasciiacute}=d\right)$,}

\textit{$\bullet$}\textit{\emph{if }}\emph{$d>c$,}\textit{\emph{
then}}\emph{ $P\left(X=d\right)\geq P\text{\textasciiacute}\left(X\text{\textasciiacute}=d\right)$.}

$\:$

It is easy to show that strong dominance implies weak dominance, but
that the converse is not true. Coming back to our CCP model, we propose
the following:

$\mathbf{Conjecture.}$ If $p\prec q$, then the random variable $Y_{n}$
defined on the probability space $\left(\Omega,P_{p}\right)$ strongly
stochastically dominates the random variable $Y_{n}$ defined on the
probability space $\left(\Omega,P_{q}\right).$

$\:$

This conjecture has been tested on various examples, but we have been
able to prove it formally for only a few values of the pair $\left(n,N\right)$,
namely for $n=2$ or 3 and any $N$, and for $n=4$ and $N\leq5$.

In applications, strong dominance reinforces weak dominance. It gives
more precise statements concerning the relative probabilities that
a given number of hosts are parasitized after a given number of eggs
laid, for two risk distributions. 

$\:$

\section*{Acknowledgments}

N.Z. and M.J.F.S acknowledge the financial support of the Fundación
Seneca of the Comunidad Autónoma de la Región de Murcia, project 19320/IP/14.
N.Z. is also grateful to the University François-Rabelais of Tours,
for its support and hospitality.


\begin{thebibliography}{Hernandez-Suarez and Hiebeler(2011)}
\bibitem[Anceaume et al.(2015)]{Anceaume2015}\textcolor{black}{Anceaume,
E., Busnel, Y., \& Sericola, B. (2015). New results on a generalized
coupon collector problem using Markov chains. Journal of Applied Probability,
52(2), 405-418.}

\bibitem[Anceaume et al.(2016)]{Anceaume2016}\textcolor{black}{Anceaume,
E., Busnel, Y., Shulte-Geers, E. \& Sericola, B. (2016). Optimization
results for a generalized coupon collector problem. Journal of Applied
Probability, 53(2), 622-629.}

\bibitem[Arnold(1987)]{ArnoldRobinHoodTransfer}Arnold, B. (1987).
Majorization and the Lorenz Order: A Brief Introduction. Arnold, B.
C.

\bibitem[Bailey et al.(2013)]{Baileyetal.2013}Bailey L. L., MacKenzie,
D. I., Nichols J. D. (2013). Special feature. Modelling demographic
processes in marked populations: Proceedings of the Euring 2013 Analytical
Meeting. Advances and applications of occupancy model. Methods in
Ecology and Evolution, doi:10.111111/2041-210X.12100

\bibitem[Boneh and Hofri(1989)]{BonehHofri}Boneh, A., Hofri, M. (1989).
The coupon-colector problem revisited. Computer Science Technical
Report. Paper 807. http://docs-lib.purdue.edu/cstech/807.

\bibitem[Fitzpatrick(1993)]{BungeFitzpatrick1993}Bunge, J., \& Fitzpatrick,
M. (1993). Estimating the number of species: a review. Journal of
the American Statistical Association, 88(421), 364-373.

\bibitem[Casas(1989)]{Casas1989}Casas, J. (1989). Foraging behaviour
of a leafminer parasitoid in the field. Ecological Entomology, 14(3),
257-265.

\bibitem[Casas et al(2004)]{Casasetal2004}Casas, J., Swarbrick, S.,
\& Murdoch, W. W. (2004). Parasitoid behaviour: predicting field from
laboratory. Ecological entomology, 29(6), 657-665.

\bibitem[Daley et al.(2001)]{DaleyGaniGani2001}Daley, D.J., Gani,
J. \& Gani, J.M. (2001). Epidemic modelling: an introduction. Cambridge
University Press.

\bibitem[Dennehy(2009)]{Dennehy2009} Dennehy, J. J. (2009). Bacteriophages
as model organisms for virus emergence research. Trends in microbiology,
17(10), 450-457.

\bibitem[Dixon(2006)]{Dixon2006}Dixon, C. J. (2006). A means of estimating
the completeness of haplotype sampling using the Stirling probability
distribution. Molecular Ecology Notes, 6(3), 650-652.

\bibitem[Donelly(1986)]{Donelly1986}Donnelly, P. (1986). Partition
structures, Polya urns, the Ewens sampling formula, and the ages of
alleles. Theoretical Population Biology, 30(2), 271-288.

\bibitem[Doumas(2015)]{Doumas 2015}Doumas, A.V. (2015). How many
trials does it take to collect all different types of a population
with probability p?. Journal of Applied Mathematics and Bioinformatics,
5(3), 1-14.

\bibitem[Ewens(1972)]{Ewens1972}Ewens, W. J. (1972). The sampling
theory of selectively neutral alleles. Theoretical population biology,
3(1), 87-112.

\bibitem[Feller(1968)]{Feller} Feller, V. (1968). An Introduction
to Probability Theory and Its Applications: Volume One. John Wiley
\& Sons.

\bibitem[Fiske(1910)]{Fiske1910}Fiske, W. F. (1910). Superparasitism:
an important factor in the natural control of insects. Journal of
Economic Entomology, 3(1), 88-97.

\bibitem[Flajolet et al.(1992)]{FlajoletBirthdayParadox}Flajolet,
P., Gardy, D., \& Thimonier, L. (1992). Birthday paradox, coupon collectors,
caching algorithms and self-organizing search. Discrete Applied Mathematics,
39(3), 207-229.

\bibitem[Hassell(2000)]{Hassell2000}Hassell, M. (2000). The spatial
and temporal dynamics of host-parasitoid interactions. Oxford University
Press.

\bibitem[Hemerick et al.(2002)]{Hemericketal.2002}Hemerik, L., Van
der Hoeven, N., \& van Alphen, J. J. (2002). Egg distributions and
the information a solitary parasitoid has and uses for its oviposition
decisions. Acta biotheoretica, 50(3), 167-188.

\bibitem[Hernandez-Suarez and Hiebeler(2011)]{Hernandez-SuarezHiebeler2011}Hernandez-Suarez,
C. M., \& Hiebeler, D. (2012). Modeling species dispersal with occupancy
urn models. Theoretical Ecology, 5(4), 555-565.

\bibitem[Huillet and Paroissin(2009)]{HuilletParoissin2013}Huillet,
T., \& Paroissin, C. (2009). Sampling from Dirichlet partitions: estimating
the number of species. Environmetrics, 20(7), 853-876.

\bibitem[Ives et al.(1999)]{IvesSchoolerJagarKnutesonGrbicSettle1999}Ives,
A.R., Schooler, S.S., Jagar, V.J., Grbic, M. \&Settle, W.H. (1999).
Variability and parasitoid foraging efficiency: a case study of pea
aphids and Aphidius ervi. The American Naturalist, 154(6), 652-673.

\bibitem[Kershenbaum et al.(2015)]{Kershenbaumetal.2015}Kershenbaum
A, Freeberg TM, Gammon DE. 2015. Estimating vocal repertoire size
is like collecting coupons: A theoretical framework with heterogeneity
in signal abundance. Journal of Theoretical Biology, 373(2015): 1-11.

\bibitem[Keeling and Rohani(2008)]{keelingRohani2008}Keeling, M.J.,
\& Rohani, P. (2008). Modelling infectious diseases in humans and
animals. Ptinceton University Press.

\bibitem[Lloyd-Smith et al.(2005)]{Lloyd-SmithSchreiberKoppGetz(2005)}
Lloyd-Smith, J.O., Schreiber, S.J., Kopp, P.E. \& Getz, W.M. (2005).
Superspreading and the effect of individual variation on disease emergence.
Nature, 438(7066), 355-359.

\bibitem[Marshall et al.(2011)]{Majorization}\textcolor{black}{Marshall,
A. W.,Olkin, I., Arnold, B. C. (2011) Inequalities: Theory of majorization
and its applications. Second edition. Springer Series in Statistics.
Springer, New York.}

\bibitem[May(1978)]{May1978}May, R.M. (1978). Host-parasitoid systems
in patchy environments: a phenomenological model. The Journal of Animal
Ecology, 833-844.

\bibitem[McArthur(1957)]{McArthur}MacArthur, R. H. (1957). On the
relative abundance of bird species. Proceedings of the National Academy
of Sciences of the United States of America, 43(3), 293-295.

\bibitem[Montovan et al.(2015) ]{Montovan} Montovan, K.J., Couchoux,
C., Jones, L.E., Reeve, H.K., van Nouhuys, S. (2015) The puzzle of
partial resource use by a parasitoid wasp. The American Naturalist,
185(4) 538-550.

\bibitem[Muirhead(1902)]{Muirhead1903}Muirhead, R. F. (1902). Some
methods applicable to identities and inequalities of symmetric algebraic
functions of n letters. Proceedings of the Edinburgh Mathematical
Society, 21, 144-162.

\bibitem[Murdoch et al.(2013)]{MurdochBriggsNisbet2013}Murdoch, W.W.,
Briggs, C.J. \& Nisbet, R.M. (2013). Consumer-Resource Dynamics (MPB-36).
Princeton University Press.

\bibitem[Neal and Moriary(2009)]{NealMoriary2009}Neal, P., \& Moriary,
J. (2009). Sampling Efficiency and Biodiversity.

\bibitem[Simpson(1949)]{Simpson}Simpson, E. H. (1949). Measurement
of diversity. Nature, 163, 688.

\bibitem[Singh et al.(2009)]{SinghMurdochNisbet2009}Singh, A., Murdoch,
W.W. \& Nisbet, R.M. (2009). Skewed attacks, stability, and host suppression.
Ecology, 90(6), 1679-1686.

\bibitem[Tenaillon et al.(2012)]{Tenaillonetal.2012}Tenaillon, O.,
Rodríguez-Verdugo, A., Gaut, R. L., McDonald, P., Bennett, A. F.,
Long, A. D., \& Gaut, B. S. (2012). The molecular diversity of adaptive
convergence. Science, 335(6067), 457-461.

\bibitem[Thompson(1924)]{Thompson1924}Thompson, W. R. (1924). La
théorie mathématique de l'action des parasites entomophages et le
facteur du hasard. Ann. Fac. Sci. Marseille, 2(2), 69-89.

\bibitem[ Vandewalle et al.(2015)]{Vandewalleetal.2015}Vandewalle,
K., Festjens, N., Plets, E., Vuylsteke, M., Saeys, Y., \& Callewaert,
N. (2015). Characterization of genome-wide ordered sequence-tagged
Mycobacterium mutant libraries by Cartesian Pooling-Coordinate Sequencing.
Nature communications, 6.

\bibitem[Wong and Yue(1973)]{WongYue1973WeakDominance}Wong, C. K.,
\& Yue, P. C. (1973). A majorization theorem for the number of distinct
outcomes in N independent trials. Discrete Mathematics, 6(4), 391-398.\end{thebibliography}
\end{document}